# GENERAL EMPIRICAL BAYES WAVELET METHODS AND EXACTLY ADAPTIVE MINIMAX ESTIMATION[1]

By Cun-Hui Zhang

*Rutgers University*

In many statistical problems, stochastic signals can be represented as a sequence of noisy wavelet coefficients. In this paper, we develop general empirical Bayes methods for the estimation of true signal. Our estimators approximate certain oracle separable rules and achieve adaptation to ideal risks and exact minimax risks in broad collections of classes of signals. In particular, our estimators are uniformly adaptive to the minimum risk of separable estimators and the exact minimax risks simultaneously in Besov balls of all smoothness and shape indices, and they are uniformly superefficient in convergence rates in all compact sets in Besov spaces with a finite secondary shape parameter. Furthermore, in classes nested between Besov balls of the same smoothness index, our estimators dominate threshold and James–Stein estimators within an infinitesimal fraction of the minimax risks. More general block empirical Bayes estimators are developed. Both white noise with drift and nonparametric regression are considered.

**1. Introduction.** Suppose a sequence $y \equiv \{y_{jk}\}$ of infinite length is observed, with

(1.1) $\quad y_{jk} \equiv \beta_{jk} + \varepsilon z_{jk}, \qquad 1 \le k \le \max(2^j, 1), \ j = -1, 0, 1, \ldots,$

where $\varepsilon > 0$ and $z_{jk}$ are i.i.d. $N(0,1)$. In many statistical problems, stochastic signals can be represented in the form of (1.1) as noisy wavelet coefficients with errors $\varepsilon z_{jk}$, or simply represented by a sequence of normal variables as in (2.2) below. In this paper we consider estimation of the true wavelet

Received August 2000; revised November 2002.

[1]Supported in part by the National Science Foundation and the National Security Agency.

coefficients $\beta \equiv \{\beta_{jk}\}$, that is, the normal means, with the $\ell_2$ risk

$$(1.2) \qquad R^{(\varepsilon)}(\hat{\beta}, \beta) \equiv \sum_{j=-1}^{\infty} \sum_{k=1}^{2^j \vee 1} E_\beta^{(\varepsilon)}(\hat{\beta}_{jk} - \beta_{jk})^2$$

for estimates $\hat{\beta} \equiv \{\hat{\beta}_{jk}\}$ based on $y$, where $E_\beta^{(\varepsilon)}$ is the expectation in model (1.1).

We develop general empirical Bayes (GEB) estimators $\hat{\beta}^{(\varepsilon)} \equiv \hat{\beta}^{(\varepsilon)}(y)$, defined in Section 2, such that under certain mild conditions on the sequence $\beta$ the risks of $\hat{\beta}^{(\varepsilon)}$ satisfy

$$(1.3) \qquad R^{(\varepsilon)}(\hat{\beta}^{(\varepsilon)}, \beta) \approx R^{(\varepsilon,*)}(\beta) \equiv \inf_{\hat{\beta} \in \mathcal{D}^s} R^{(\varepsilon)}(\hat{\beta}, \beta),$$

where $\mathcal{D}^s$ is the class of all *separable estimators* of the form $\hat{\beta}_{jk} = h_j(y_{jk})$ with Borel $h_j$. We provide *oracle inequalities*, that is, upper bounds for the *regret* $R^{(\varepsilon)}(\hat{\beta}^{(\varepsilon)}, \beta) - R^{(\varepsilon,*)}(\beta)$ for this adaptation to the *ideal risk* $R^{(\varepsilon,*)}(\beta)$. Our oracle inequalities imply that the ideal adaptation (1.3) is uniform for large collections $\mathcal{B}$ of classes $B$ of the unknown $\beta$, for example, Lipschitz, Sobolev and Besov balls $B$ of all smoothness and shape indices and radii, in the sense that for all $B \in \mathcal{B}$ the regret is uniformly of smaller order than the minimax risk

$$(1.4) \qquad \mathcal{R}^{(\varepsilon)}(B) \equiv \inf_{\hat{\beta}} \sup_{\beta \in B} R^{(\varepsilon)}(\hat{\beta}, \beta).$$

This uniform ideal adaptation implies: (1) the *exact minimax adaptation*

$$(1.5) \qquad \sup\{R^{(\varepsilon)}(\hat{\beta}^{(\varepsilon)}, \beta) : \beta \in B\} = (1 + o(1))\mathcal{R}^{(\varepsilon)}(B)$$

simultaneously for all Besov balls $B$, (2) adaptation to spatial inhomogeneity of the signal [Donoho and Johnstone (1994a)], (3) the superefficiency of the GEB estimators in convergence rates in all compact sets of $\beta$ in Besov spaces with a finite secondary shape parameter and (4) dominance of GEB estimators over other empirical Bayes (EB) or separable estimators in the limit in all classes of $\beta$ nested between Besov balls of the same smoothness index. We also describe more general block EB methods and implementation of GEB estimators in nonparametric regression models with possibly unknown variance.

The white noise model (1.1) is a wavelet representation of its original form [cf. Ibragimov and Khas'minskii (1981)], in which one observes

$$(1.6) \qquad Y(t) \equiv \int_0^t f(u)\,du + \varepsilon W(t), \qquad 0 \le t \le 1,$$

where $f \in L^2[0,1]$ is unknown and $W(\cdot)$ is a standard Brownian motion. In this representation, $y_{jk} \equiv \int \phi_{jk}(t) \, dY(t)$, $\beta_{jk} \equiv \beta_{jk}(f) \equiv \int f(t) \phi_{jk}(t) \, dt$, and estimates

$$(1.7) \qquad \hat{f}(t) \equiv \sum_{j,k} \hat{\beta}_{jk} \phi_{jk}(t), \qquad 0 \leq t \leq 1,$$

are constructed based on estimates $\hat{\beta}_{jk}$ of $\beta_{jk}$, where $\phi_{jk}$ are wavelets forming an orthonormal basis in $L^2[0,1]$. Let $E_f^{(\varepsilon)}$ be the expectation in model (1.6). By the Parseval identity,

$$(1.8) \quad R^{(\varepsilon)}(\hat{\beta}, \beta) = E_\beta^{(\varepsilon)} \sum_{j,k} (\hat{\beta}_{jk} - \beta_{jk})^2 = E_f^{(\varepsilon)} \int_0^1 \{\hat{f}(t) - f(t)\}^2 \, dt,$$

so that our problem is equivalent to the estimation of $f$ under the mean integrated squared error (MISE). In general, $\phi_{jk}(t)$ are "periodic" or "boundary adjusted" dilation and translation $2^{j/2} \phi(2^j t - k)$ of a "mother wavelet" $\phi$ of compact support for $j \geq j_0$ for certain $j_0 \geq 0$; see Donoho and Johnstone (1994a). Since $\phi_{jk}$ is supported in an interval of size $O(1/2^j)$ in the vicinity of $k/2^j$, $j$ and $k$ are, respectively, resolution and spatial indices, and $y_{jk}$ represent the information about the behavior of $f$ at resolution level $j$ and location $k/2^j$. We refer to Chui (1992), Daubechies (1992) and Härdle, Kerkyacharian, Picard and Tsybakov (1998) for wavelet theory and its applications.

This paper is organized as follows. We develop block EB methods in Section 2 which naturally lead to GEB estimators. We state main properties of the GEB estimators in Section 3. We implement GEB estimators in nonparametric regression models in Section 4. We discuss related results and problems in Section 5. We focus on compound estimation of normal means in Section 6. We present our main theorems in their full strength in Section 7. We cover Bayes models and more general classes of the unknown $\beta$ in Section 8. We study the equivalence between the nonparametric regression and white noise models in Section 9. Proofs are given in the Appendix unless otherwise stated or provided immediately after the statements of results. The main theorems in Sections 3, 6 and 7 have been reported earlier in Zhang (2000) with more details in proofs. We use the notation $\log_+ x \equiv 1 \vee \log x$ and $x_{(n)} \equiv (x_1, \ldots, x_n)$ throughout.

**2. Block EB methods and GEB estimators.** We begin with block EB methodologies, which naturally lead to GEB estimators. Consider a sequence of $N \leq \infty$ decision problems with observations $X_k \sim p(x|\theta_k)$ and parameters $\theta_k$ under the compound risk $\sum_{k=1}^N EL_0(\delta_k, \theta_k)$ for a given loss $L_0(\cdot, \cdot)$. Block EB methods partition the sequence into blocks $[j] \equiv (k_{j-1}, k_j]$, $k_{j-1} < k_j <$ $\infty$, and apply EB procedures of the form $\delta_k = \hat{t}_{[j]}(X_k)$, $k \in [j]$, in individual blocks, where $\hat{t}_{[j]}(\cdot)$ are estimates of the *oracle rules*

$$(2.1) \qquad t^*_{[j]}(\cdot) \equiv \arg\min_{t(\cdot) \in \mathcal{D}_0} \sum_{k \in [j]} E_{\theta_k} L_0(t(X_k), \theta_k)$$

for a certain class $\mathcal{D}_0$ of decision rules. Block GEB (linear, threshold EB) procedures approximate the oracle rules (2.1) corresponding to the classes $\mathcal{D}_0$ of all Borel (linear, threshold) functions. It follows from compound decision theory [Robbins (1951)] that $t^*_{[j]}$ are the Bayes rules when the priors are taken to be the unknown empirical distributions of $\{\theta_k, k \in [j]\}$.

Consider the estimation of normal means $\beta_k$ based on independent observations

$$(2.2) \qquad y \equiv \{y_k, k \leq N\}, \qquad y_k \sim N(\beta_k, \varepsilon^2)$$

with known $\varepsilon$, under the squared loss as in (1.2). After standardization with $(X_k, \theta_k) \equiv (y_k, \beta_k)/\varepsilon$ to the unit variance, block GEB estimators of $\beta_k$ become

$$(2.3) \qquad \hat{\beta}^{(\varepsilon)}_k \equiv \hat{\beta}^{(\varepsilon)}_k(y) \equiv \varepsilon \hat{t}_{[j]}(y_k/\varepsilon), \qquad k \in [j],$$

where $\hat{t}_{[j]}$ are estimates of (2.1) with squared-error loss $L_0(\delta, \theta) = (\delta - \theta)^2$. The empirical distribution of $\{\theta_k, k \in [j]\}$ is

$$(2.4) \qquad G^{(\varepsilon)}_{[j]}(u) \equiv \frac{1}{n_j} \sum_{k \in [j]} I\{\theta_k \leq u\}, \qquad \theta_k \equiv \frac{\beta_k}{\varepsilon},$$

where $n_j \equiv k_j - k_{j-1}$ is the size of block $j$. Let $\varphi(x) \equiv e^{-x^2/2}/\sqrt{2\pi}$ and

$$(2.5) \qquad \varphi_G(x) \equiv \varphi(x; G) \equiv \int \varphi(x-u)\, dG(u), \qquad \varphi'(x; G) \equiv \frac{d\varphi_G}{dx}.$$

The oracle rules (2.1) are explicit functionals of the mixture marginal distributions $\varphi(x; G^{(\varepsilon)}_{[j]})$ of the observations $\{X_k \equiv y_k/\varepsilon, k \in [j]\}$ [Robbins (1956), page 162, and Stein (1981)], given by

$$(2.6) \qquad t^*_{[j]}(x) = x + \frac{\varphi'(x; G^{(\varepsilon)}_{[j]})}{\varphi(x; G^{(\varepsilon)}_{[j]})}.$$

This formula motivated the GEB estimators of Zhang (1997).

We construct GEB estimators in individual blocks using a hybrid version of the GEB estimator of Zhang (1997). The hybrid GEB estimator utilizes an estimate of the order of $\kappa(G^{(\varepsilon)}_{[j]})$,

$$(2.7) \qquad \kappa(G) \equiv \int (|u|^2 \wedge 1)\, dG(u),$$

and switches from the GEB estimator to a threshold estimator when $\kappa(G_{[j]}^{(\varepsilon)})$ is small. Specifically, for certain $\rho(n) > 0$ and $b(n)$ given in (2.11) below and $n_j \geq n_* > 2$, we define (2.3) by

$$(2.8) \quad \hat{t}_{[j]}(x) \equiv \begin{cases} x + \hat{\varphi}'_{[j]}(x)/\max\{\rho(n_j), \hat{\varphi}_{[j]}(x)\}, & \text{if } \hat{\kappa}_{[j]} > b(n_j), \\ \operatorname{sgn}(x)(|x| - \sqrt{2\log n_j})^+, & \text{if } \hat{\kappa}_{[j]} \leq b(n_j), \end{cases}$$

where $\hat{\varphi}_{[j]}(x)$, a kernel estimate of $\varphi(x; G_{[j]}^{(\varepsilon)})$ in (2.6), is given by

$$(2.9) \qquad \hat{\varphi}_{[j]}(x) \equiv \frac{1}{n_j} \sum_{k \in [j]} \sqrt{2\log n_j} K(\sqrt{2\log n_j}(x - X_k))$$

with $K(x) \equiv \sin(x)/(\pi x)$ and $X_k \equiv y_k/\varepsilon$, and $\hat{\kappa}_{[j]}$, an estimate of the order of $\kappa(G_{[j]}^{(\varepsilon)})$, is given by

$$\hat{\kappa}_{[j]} \equiv 1 - \frac{\sqrt{2}}{n_j} \sum_{k \in [j]} \exp(-X_k^2/2).$$

For $n_j < n_*$, we choose the MLE $\hat{\beta}_k^{(\varepsilon)} \equiv y_k$ [i.e., $\hat{t}_{[j]}(x) \equiv x$ or $\rho(n) = \infty = -b(n)$ for $n < n_*$] or the James and Stein (1961) estimator for the vectors $\{\beta_k, k \in [j]\}$.

For denoising the wavelet coefficients, we identify the sequence $y \equiv \{y_{-1,1}, y_{0,1}, y_{1,1}, y_{1,2}, \ldots\}$ in (1.1) with $y \equiv \{y_k, k \leq N\}$ in (2.2) and partition them into natural blocks $[j] = (2^j, 2^{j+1}]$, $j = -1, 0, \ldots$, with a single block for each resolution level $j$. This results in

$$(2.10) \quad \hat{\beta}^{(\varepsilon)} \equiv \hat{\beta}^{(\varepsilon)}(y) \equiv \{\hat{\beta}_{jk}^{(\varepsilon)}\}, \qquad \hat{\beta}_{jk}^{(\varepsilon)} \equiv \begin{cases} y_{jk}, & \text{if } j < j_*, \\ \varepsilon \hat{t}_{[j]}(y_{jk}/\varepsilon), & \text{if } j \geq j_*, \end{cases}$$

where $j_* \equiv \max\{j : j \leq (\log n_*)/\log 2\}$ and $\hat{t}_{[j]}$ is as in (2.8) with $n_j = 2^j$,

$$\hat{\varphi}_{[j]}(x) \equiv \frac{\sqrt{2j\log 2}}{2^j} \sum_{k=1}^{2^j} K\left(\frac{x - y_{jk}/\varepsilon}{(2j\log 2)^{-1/2}}\right), \qquad K(x) \equiv \frac{\sin(x)}{\pi x},$$

and $\hat{\kappa}_{[j]} \equiv 1 - 2^{-j} \sum_{k=1}^{2^j} \sqrt{2} \exp(-(y_{jk}/\varepsilon)^2/2)$. For definiteness, we set

$$(2.11) \quad \rho \equiv \rho(n) = (1 + \eta_n)\rho_0\sqrt{2(\log n)/n}, \qquad b \equiv b(n) = b_0(\log n)/\sqrt{n},$$

with certain $\eta_n \to 0$ and positive constants $\rho_0$ and $b_0$. We simply call (2.10) GEB estimators since the blocks represent natural resolution levels in the wavelet setting.

We discuss in detail in Section 6 the construction and properties of the GEB estimators in individual blocks (resolution levels). Here we briefly describe the rationale for our choices of the "tuning parameters" for (2.8) and (2.10). The special kernel and bandwidth in (2.9) ensure that $\hat{\varphi}_{[j]}(x) \to \varphi(x; G_{[j]}^{(\varepsilon)})$ at nearly the optimal rate $n_j^{-1/2}$ uniformly and in derivatives as $n_j \to \infty$ [Zhang (1997)]. The sample size $n_*$ (and thus the initial resolution level $j_*$) should be determined so that $\hat{\varphi}_{[j]}(x) \approx \varphi(x; G_{[j]}^{(\varepsilon)})$ with sufficient accuracy for $n_j \geq n_*$. Although the $\rho(n_j)$ and $b(n_j)$ in (2.8) could be determined/optimized by data-driven methods, for example, Stein's (1981) estimator of mean squared error, bootstrap and cross validation, properties of the resulting estimators are not clear. The choice in (2.11) provides the sharpest bounds in our main theorems. Our risk bounds depend on $j_*$, $\rho_0$ and $b_0$ only through scaling constants in terms of smaller order than minimax risks. Finally, we remark that (block) GEB estimators (2.3) and (2.10) are scale equivariant:

$$(2.12) \qquad \hat{\beta}^{(\varepsilon)}(y) = C\hat{\beta}^{(\varepsilon/C)}(y/C) \qquad \forall C > 0,$$

since $\hat{t}_{[j]}(x)$ in (2.8) depend on $y$ and $\varepsilon$ only through $y/\varepsilon$. Thus, for the risks in (1.3) and all $C > 0$

$$
\begin{aligned}
R^{(\varepsilon)}(\hat{\beta}^{(\varepsilon)}(y), \beta) &= C^2 R^{(\varepsilon/C)}(\hat{\beta}^{(\varepsilon/C)}(y/C), \beta/C), \\
R^{(\varepsilon,*)}(\beta) &= C^2 R^{(\varepsilon/C,*)}(\beta/C).
\end{aligned}
\tag{2.13}
$$

**3. Oracle inequalities and their consequences.** In this section we describe main properties of our (block) GEB estimators (2.3), (2.8) and (2.10) and the concepts of uniform ideal adaptivity, exactly adaptive minimaxity, spatial adaptivity and superefficiency. Sections 5, 7 and 8 contain further discussion about these properties and concepts.

3.1. *Oracle inequalities.* Consider the estimation of normal means with observations (2.2). An *oracle expert* with the knowledge of $t_{[j]}^*$ in (2.6) could use the ideal separable rule $\varepsilon t_{[j]}^*(y_k/\varepsilon)$ for $\beta_k$ to achieve the ideal risk

$$(3.1) \quad R^{(\varepsilon,*)}(\beta) \equiv \min_{\hat{\beta} \in \mathcal{D}^s} R^{(\varepsilon)}(\hat{\beta}, \beta) = \sum_j \min_{t(\cdot)} \sum_{k \in [j]} E_\beta^{(\varepsilon)}\{\varepsilon t(y_k/\varepsilon) - \beta_k\}^2,$$

as in (1.3), where $\mathcal{D}^s$ is the collection of all separable estimates of the form $\hat{\beta}_k = h_j(y_k)$, $\forall k \in [j]$. Although $\varepsilon t_{[j]}^*(y_{jk}/\varepsilon)$ are not statistics, the ideal risk (3.1) provides a benchmark for our problem.

THEOREM 3.1. *Let $\hat{\beta}^{(\varepsilon)} \equiv \{\hat{\beta}_k^{(\varepsilon)}\}$ be as in (2.3) and (2.8) based on (2.2). Let $R^{(\varepsilon)}(\hat{\beta}, \beta) \equiv \sum_{k=1}^N E_\beta^{(\varepsilon)}(\hat{\beta}_k - \beta_k)^2$. Then there exists a universal constant*

$M < \infty$ such that

$$R^{(\varepsilon)}(\hat{\beta}^{(\varepsilon)}, \beta) - R^{(\varepsilon,*)}(\beta)$$
(3.2)
$$\leq M\varepsilon^2 \sum_j \left\{ n_j r_{p\wedge 2}\left(n_j, \frac{(\sum_{k\in[j]}|\beta_k|^p)^{1/p}}{\varepsilon n_j^{1/p}}\right) + \frac{1}{(\log n_j + 1)^{3/2}} \right\},$$

where $R^{(\varepsilon,*)}(\beta)$ is the ideal risk in (3.1), $n_j \equiv k_j - k_{j-1}$ are block sizes and

$$(3.3)\ r_p(n, C) \equiv \min\left(1, \frac{C^p}{(\log n)^{p/2-1}}, \max\left[\frac{(\log n)^2}{\sqrt{n}}, \left\{\frac{C(\log n)^{3/2}}{\sqrt{n}}\right\}^{p/(p+1)}\right]\right).$$

COROLLARY 3.1. *If* $\beta = 0$, *then* $R^{(\varepsilon,*)}(\beta) = 0$ *and* $R^{(\varepsilon)}(\hat{\beta}^{(\varepsilon)}, \beta) = O(\varepsilon^2)$.

Theorem 3.1, proved in Section 7.1, provides a crucial oracle inequality in the derivation of our main results. It allows us to bound the regret of our estimators in terms of the moments of $\beta_{[j]}$. Consider block sizes $n_j$ such that, for all $p > 0$ and $\eta > 0$ and as $x \to \infty$,

$$(3.4)\quad \sum_{n_j \geq x} \frac{x^p}{n_j^p} = o(x^\eta), \qquad \sum_{n_j \leq x} \frac{n_j^p}{x^p} = o(x^\eta), \qquad \sum_j (1 + \log n_j)^{-3/2} < \infty.$$

Condition (3.4) holds if $\log n_j \sim j^\gamma$ for certain $2/3 < \gamma \leq 1$.

THEOREM 3.2. *Let* $\hat{\beta}^{(\varepsilon)}$ *be as in Theorem* 3.1 *and* $\|\beta\| \equiv \sup_j n_j^s \times (\sum_{k\in[j]}|\beta_k|^p/n_j)^{1/p}$. *Suppose* (3.4) *holds and* $\alpha_0 \equiv \min\{s, 2s - 1/2 - 1/(p \wedge 2)\} > 0$. *Then, for all* $\eta > 0$

$$\sup\{R^{(\varepsilon)}(\hat{\beta}^{(\varepsilon)}, \beta) - R^{(\varepsilon,*)}(\beta) : \|\beta\| \leq C\} \leq o(\varepsilon^{2\alpha_0/(\alpha_0+1/2)-\eta}) \qquad as\ \varepsilon \to 0+.$$

REMARK. If a higher threshold level $\sqrt{2(1+A_0)\log n_j}$ with $A_0 > 0$ is used in (2.8) for $\hat{\kappa}_{[j]} \leq b(n_j)$, Theorems 3.1 and 3.2 hold with $(1+\log n_j)^{-3/2}$ replaced by $n_j^{-A_0}(1+\log n_j)^{-3/2}$ in (3.2) and (3.4). See the remark below Theorem 6.4.

In the rest of Section 3, we focus on the wavelet model (1.1), that is, the case of $n_j = 2^j$. Our methodology is clearly applicable to more general block sizes $n_j$ satisfying (3.4).

3.2. *Uniform ideal adaptation.* Let $R^{(\varepsilon)}(\hat{\beta}, \beta)$ be the $\ell_2$ risk in (1.2). Statistical estimators $\hat{\beta}^{(\varepsilon)}$ are uniformly adaptive to the ideal risk $R^{(\varepsilon,*)}(\beta)$ in (1.3) and (3.1), with respect to a collection $\mathcal{B}$ of classes $B$ of the unknown sequence $\beta$, if

$$(3.5)\ \sup_{\beta \in B}\{R^{(\varepsilon)}(\hat{\beta}^{(\varepsilon)}, \beta) - R^{(\varepsilon,*)}(\beta)\} = o(1)\mathcal{R}^{(\varepsilon)}(B) \qquad as\ \varepsilon \to 0+\ \forall\, B \in \mathcal{B},$$

and $\hat{\beta}^{(\varepsilon)}$ depends on $(y,\varepsilon)$ only, not on $B$, where $\mathcal{R}^{(\varepsilon)}(B)$ is the minimax risk in (1.4). In other words, uniform ideal adaptation demands that, for all $B \in \mathcal{B}$ and in the minimax sense, the *regret*

$$(3.6) \qquad r^{(\varepsilon)}(\hat{\beta}^{(\varepsilon)}, \beta) \equiv R^{(\varepsilon)}(\hat{\beta}^{(\varepsilon)}, \beta) - R^{(\varepsilon,*)}(\beta)$$

be uniformly of smaller order than the typical convergence rates in $B$. As an immediate consequence of uniform ideal adaptation, maximum risks are bounded by the maximum ideal risks,

$$(3.7) \qquad \sup_{\beta \in B} R^{(\varepsilon)}(\hat{\beta}^{(\varepsilon)}, \beta) \leq (1 + o(1)) \sup_{\beta \in B} R^{(\varepsilon,*)}(\beta) \qquad \forall\, B \in \mathcal{B}.$$

Our GEB estimators possess this uniform ideal adaptivity property with respect to

$$(3.8) \qquad \mathcal{B}_{\text{Besov}} \equiv \bigg\{ B_{p,q}^{\alpha}(C) : 0 < \alpha < \infty,$$

$$\frac{1}{\alpha + 1/2} < p \leq \infty, 0 < q \leq \infty, 0 < C < \infty \bigg\},$$

where $B_{p,q}^{\alpha} \equiv B_{p,q}^{\alpha}(C)$ are the Besov balls defined by

$$(3.9) \qquad B_{p,q}^{\alpha} \equiv \{\beta : \|\beta\|_{p,q}^{\alpha} \leq C\},$$

$$\|\beta\|_{p,q}^{\alpha} \equiv \left[ |\beta_{-1,0}|^q + \sum_{j=0}^{\infty} (2^{j(\alpha+1/2-1/p)} \|\beta_{[j]}\|_{p,2^j})^q \right]^{1/q},$$

with $\|\beta_{[j]}\|_{p,2^j} \equiv (\sum_{k=1}^{2^j} |\beta_{jk}|^p)^{1/p}$, and with the usual modifications for $p \vee q = \infty$. For $p \wedge q < 1$, $\|\cdot\|_{p,q}^{\alpha}$ is not a norm, but $(\|\beta' + \beta''\|_{p,q}^{\alpha})^{p \wedge q} \leq (\|\beta'\|_{p,q}^{\alpha})^{p \wedge q} + (\|\beta''\|_{p,q}^{\alpha})^{p \wedge q}$ is sufficient here.

THEOREM 3.3. *Let $\hat{\beta}^{(\varepsilon)} \equiv \{\hat{\beta}_{jk}^{(\varepsilon)}\}$ be as in (2.10) based on $y \equiv \{y_{jk}\}$ in (1.1), with $\rho(n)$ and $b(n)$ in (2.11). Then (3.5) holds for $\mathcal{B} = \mathcal{B}_{\text{Besov}}$.*

By Donoho and Johnstone [(1998), Theorem 1] and Theorem 7.3 below, the minimax convergence rates in Besov balls are given by

$$(3.10) \ 0 < \inf_{0 < \varepsilon \leq C} \frac{\mathcal{R}^{(\varepsilon)}(B_{p,q}^{\alpha}(C))}{\varepsilon^{2\alpha/(\alpha+1/2)} C^{1/(\alpha+1/2)}} \leq \sup_{0 < \varepsilon \leq C} \frac{\mathcal{R}^{(\varepsilon)}(B_{p,q}^{\alpha}(C))}{\varepsilon^{2\alpha/(\alpha+1/2)} C^{1/(\alpha+1/2)}} < \infty.$$

Based on (3.10), Theorem 3.3 is an immediate consequence of Theorems 3.2 and 7.2 in Section 7, which provide upper bounds for the convergence rates of the $o(1)$ in (3.5). Note that $\alpha_0 > \alpha \equiv s - 1/2$ in Theorem 3.2 for $s > 1/p$. We show in Section 8 that (3.5) and (1.5) hold for much larger collections than $\mathcal{B}_{\text{Besov}}$.

3.3. *Adaptive minimaxity.* A main consequence of the uniform ideal adaptivity in Theorem 3.3 is the universal exactly adaptive minimaxity over all Besov balls.

THEOREM 3.4. *Let $\hat{\beta}^{(\varepsilon)} \equiv \{\hat{\beta}^{(\varepsilon)}_{jk}\}$ be as in (2.10) and (2.11) with positive constants $(j_*, \rho_0, b_0)$. Then (1.5) holds for the Besov balls $B^\alpha_{p,q}(C)$ in (3.9) for all $(\alpha, p, q, C)$ in (3.8).*

This result can be viewed as an extension of the work of Efromovich and Pinsker (1984, 1986), Efromovich (1985) and Golubev (1992) from Sobolev versions of $B_{\alpha,2,2}$ to Besov balls with general $(\alpha, p, q)$. Theorem 3.4 follows from Theorem 7.4 in Section 7, which provides upper bounds for the order of the $o(1)$ in (1.5).

For a general collection $\mathcal{B}$, exact adaptive minimaxity (1.5) is a consequence of (3.5) and

(3.11) $$\sup_{\beta \in B} R^{(\varepsilon,*)}(\beta) = (1 + o(1))\mathcal{R}^{(\varepsilon)}(B),$$

since (3.5) implies (3.7). For Besov balls $B = B^\alpha_{p,q}$, (3.11) is proved in Donoho and Johnstone (1998) for $q \geq p$ and in Theorem 7.3 for general $(p, q)$.

3.4. *Spatial adaptation.* Another main consequence of the uniform ideal adaptivity in Theorem 3.3 is spatial adaptivity of (1.7) when $\beta \equiv \beta(f)$ represents wavelet coefficients of a spatially inhomogeneous signal function $f(\cdot)$. For $\beta \in B^\alpha_{p,q}$, the *smoothness index* $\alpha$ indicates the typical rate of decay of $|\beta_{jk}|$ as $j \to \infty$. Donoho and Johnstone (1994a) and Donoho, Johnstone, Kerkyacharian and Picard (1995) pointed out that spatial inhomogeneity of a function $f$ is often reflected in the sparsity of its wavelet coefficients $\beta_{jk} \equiv \beta_{jk}(f)$ at individual resolution levels, not necessarily in the smoothness index $\alpha$. In such cases, a handful of $|\beta_{jk}|$ could be much larger than the overall order of magnitude of $\beta_{[j]}$ at individual resolution levels, so that $\beta \in B^\alpha_{p,q}$ only for small $p < 2$. Thus, spatial adaptation can be achieved via (exactly, *rate* or *nearly*) adaptive minimaxity in Besov balls with small shape parameter $p$. Our GEB estimators are spatially adaptive to the full extent in the sense that they are exactly adaptive minimax in Besov balls for all $(\alpha, p, q)$, under the minimum condition $p > 1/(\alpha + 1/2)$, even allowing $p < 1$.

EXAMPLE 3.1. Let $\mathcal{F}_{d,m}(C)$ be the collection of all piecewise polynomials $f$ of degree $d$ in $[0, 1]$, with at most $m$ pieces and $\|f\|_\infty \leq C$. Let $\phi$ be a mother wavelet with $\int x^j \phi(x)\,dx = 0$, $j = 0, \ldots, d$, and $\phi(x) = 0$ outside an interval $I_0$ of length $|I_0|$. For $f \in \mathcal{F}_{d,m}(C)$, the wavelet coefficients $\beta_{jk} \equiv \beta_{jk}(f) \equiv 2^{j/2} \int_0^1 f(x)\phi(2^j x - k)\,dx = 0$ if $f$ is a single piece of polynomial in $(I_0 + k)/2^j$ and $|\beta_{jk}| \leq 2^{-j/2} C \int |\phi|\,dx$ otherwise. Thus, $\|\beta_{[j]}\|_{p,2^j} \leq 2^{-j/2} m^{1/p} C M_0$ for all $j$ and $p$, where $M_0 \equiv (|I_0| + 2)^{1/p} \int |\phi|\, dx$. By (3.9), $\|\beta\|_{p,q}^\alpha < \infty$ if $\alpha < 1/p$ for $q < \infty$ or $\alpha = 1/p$ for $q = \infty$. Theorem 3.4, (3.10) and (1.8) imply that $E_f^{(\varepsilon)} \int_0^1 (\hat{f} - f)^2\, dx = R^{(\varepsilon)}(\hat{\beta}^{(\varepsilon)}, \beta(f)) = O(\varepsilon^{2\alpha/(\alpha+1/2)})$ for all $\alpha < \infty$. Moreover, Theorem 8.1 in Section 8 implies that for $m^{(\varepsilon)} = o(1)(\log \varepsilon)^{-2} \varepsilon^{-1/(\alpha+1/2)}$,

$$\limsup_{\varepsilon \to 0+} \varepsilon^{-2\alpha/(\alpha+1/2)} \sup \{R^{(\varepsilon)}(\hat{\beta}^{(\varepsilon)}, \beta(f)) : f \in \mathcal{F}_{d,m^{(\varepsilon)}}(\varepsilon^{-M})\} = 0 \qquad \forall M < \infty,$$

with the radii $C = 0$ in (8.4).

3.5. *Superefficiency.* An interesting phenomenon with our GEB estimators is their universal *superefficiency* in convergence rates in compact sets in Besov spaces with $q < \infty$.

THEOREM 3.5. *Let $\hat{\beta}^{(\varepsilon)} \equiv \{\hat{\beta}_{jk}^{(\varepsilon)}\}$ be as in (2.10) and (2.11) with positive constants $(j_*, \rho_0, b_0)$. Let $0 < \alpha < \infty$, $1/(\alpha+1/2) < p \leq \infty$ and $0 < q < \infty$. Then $\lim_{\varepsilon \to 0+} \varepsilon^{-2\alpha/(\alpha+1/2)} R^{(\varepsilon)}(\hat{\beta}^{(\varepsilon)}, \beta) = 0$ for $\|\beta\|_{p,q}^\alpha < \infty$, and for $\|\cdot\|_{p,q}^\alpha$-compact sets $B$*

(3.12) $$\lim_{\varepsilon \to 0+} \varepsilon^{-2\alpha/(\alpha+1/2)} \sup \{R^{(\varepsilon)}(\hat{\beta}^{(\varepsilon)}, \beta) : \beta \in B\} = 0.$$

Theorem 3.5 is proved at the end of Section 7. It indicates that the minimax risks $\mathcal{R}^{(\varepsilon)}(B_{p,q}^\alpha) \sim \varepsilon^{2\alpha/(\alpha+1/2)}$ are quite conservative as measurements of the risk of our GEB estimators. As a function of $\beta$, the ideal risk $R^{(\varepsilon,*)}(\beta)$ provides more accurate information about the actual risk; see Theorems 3.2 and 7.2. Brown, Low and Zhao (1997) constructed universal pointwise superefficient estimators for Sobolev spaces (i.e., $p = 2$). Their method also provides the superefficiency of the estimators of Efromovich and Pinsker (1984, 1986). The classical kernel and many other smoothing methods do not possess the superefficiency property. In parametric models superefficiency could possibly happen only in a very small part of the parameter space, while the superefficiency of the GEB estimators is universal in all Besov balls.

3.6. *Dominance of GEB methods.* Consider classes $B$ of the unknown $\beta$ satisfying

(3.13) $$B \subseteq B_{p,q}^\alpha(C), \qquad \liminf_{\varepsilon \to 0+} \varepsilon^{-2\alpha/(\alpha+1/2)} \mathcal{R}^{(\varepsilon)}(B) > 0,$$

for certain $(\alpha, p, q, C)$ in (3.8), where $\mathcal{R}^{(\varepsilon)}(B)$ is the minimax risk (1.4). It follows from Theorem 3.4 that our GEB estimators achieve the minimax rate of convergence in $B$, but they may not achieve the minimax constant for $B$

GENERAL EB METHODS 11

in the limit. We show here that the GEB estimators dominate restricted EB estimators within $o(1)\varepsilon^{2\alpha/(\alpha+1/2)}$ in risk in all classes $B$ satisfying (3.13).

Let $\widetilde{R}^{(\varepsilon,*)}(\beta)$ be certain "ideal risk" with $\widetilde{R}^{(\varepsilon,*)}(\beta) \geq R^{(\varepsilon,*)}(\beta)$ and consider $\widetilde{\beta}^{(\varepsilon)}$ satisfying

$$(3.14) \qquad \sup_{\beta \in B} \{\widetilde{R}^{(\varepsilon,*)}(\beta) - R^{(\varepsilon)}(\tilde{\beta}^{(\varepsilon)}, \beta)\} \leq o(1)\varepsilon^{2\alpha/(\alpha+1/2)}.$$

THEOREM 3.6. *Let $\hat{\beta}^{(\varepsilon)} \equiv \{\hat{\beta}_{jk}^{(\varepsilon)}\}$ be as in* (2.10) *and* (2.11) *with positive constants* $(j_*, \rho_0, b_0)$. *Let $\mathcal{R}^{(\varepsilon)}(B)$ be the minimax risk in* (1.4). *Suppose* (3.13) *and* (3.14) *hold. Then*

$$(3.15) \qquad \lim_{\varepsilon \to 0+} \frac{\sup\{R^{(\varepsilon)}(\hat{\beta}^{(\varepsilon)}, \beta) - R^{(\varepsilon)}(\tilde{\beta}^{(\varepsilon)}, \beta) : \beta \in B\}}{\mathcal{R}^{(\varepsilon)}(B)} \leq 0.$$

*Consequently,* $\lim_{\varepsilon \to 0+} \{\sup_{\beta \in B} R^{(\varepsilon)}(\hat{\beta}^{(\varepsilon)}, \beta)\} / \{\sup_{\beta \in B} R^{(\varepsilon)}(\tilde{\beta}^{(\varepsilon)}, \beta)\} \leq 1$.

Theorem 3.6 is an immediate consequence of Theorem 3.3. Condition (3.13) holds if $B = \{\beta : \|\beta\| \leq C\}$ are balls for a certain norm $\|\cdot\|$ nested between two Besov norms with $M^{-1}\|\beta\|_{p',q'}^\alpha \leq \|\beta\| \leq M\|\beta\|_{p,q}^\alpha$ for a certain $0 < M < \infty$, for example, Lipschitz and Sobolev classes. Examples of $\tilde{\beta}^{(\varepsilon)}$ satisfying (3.14) include the Johnstone and Silverman (1998, 2005) parametric EB, block threshold (e.g., VISUAL- and SURESHRINK) and linear (e.g., James–Stein) estimators with

$$(3.16) \qquad \widetilde{R}^{(\varepsilon,*)}(\beta) \equiv \sum_j \inf_{t \in \mathcal{D}_0} \sum_k E_\beta^{(\varepsilon)}(t(y_{jk}) - \beta_{jk})^2$$

for restricted classes $\mathcal{D}_0$ (e.g., threshold, linear) of functions $t(\cdot)$.

**4. Nonparametric regression.** In this section we describe implementation of our GEB estimators in the nonparametric regression model

$$(4.1) \qquad Y_i \equiv f(t_i) + e_i, \qquad e_i \sim N(0, \sigma^2),\ i \leq N.$$

We report some simulation results for $t_i = i/N$ and unknown variance $\sigma^2$, and present the exact adaptive minimaxity and superefficiency of GEB estimators for i.i.d. uniform $t_i$ and known $\sigma^2$.

4.1. *Deterministic design and simulation results.* The white noise model (1.1) is directly connected to the nonparametric regression model via discrete wavelet reconstruction. Suppose $t_i = i/N$ and $N = 2^{J+1}$ in (4.1). A discrete wavelet reconstruction can be expressed by invertible linear mappings

$$(4.2) \qquad (y_{jk}, k \leq 2^{j \vee 0}, j \leq J) = N^{-1/2}\mathcal{W}_{N \times N}(Y_i, i \leq N),$$

$$Y_i = \sqrt{N} \sum_{j=-1}^{J} \sum_{k=1}^{2^{j \vee 1}} y_{jk} W_{jk}(i),$$

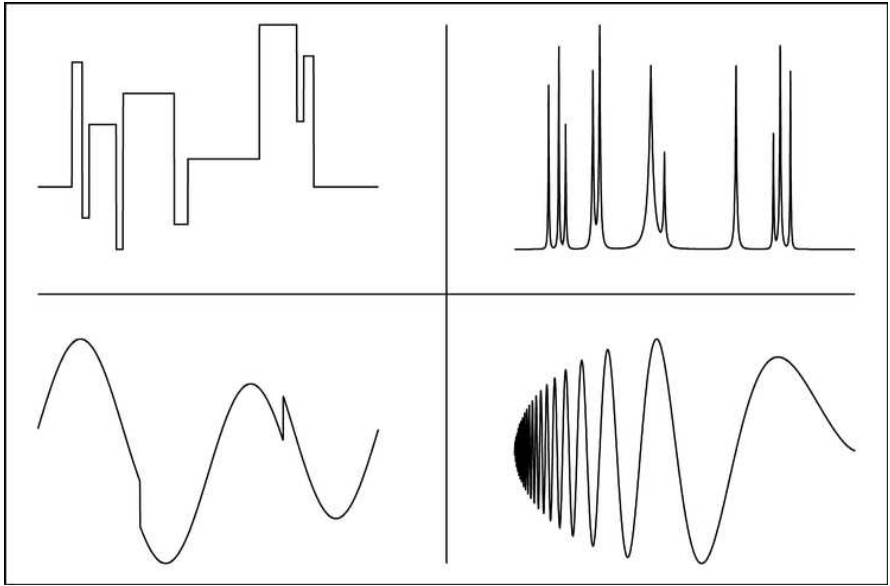

Fig. 1. *Signals; clockwise from top left, Blocks, Bumps, HeaviSine, Doppler.*

where $\mathcal{W}_{N\times N}$, called the finite wavelet transformation matrix, is a real orthonormal matrix, $W_{jk}(i)$ specify the inverse of $\mathcal{W}_{N\times N}$, and $\sqrt{N}W_{jk}(i) \approx \phi_{jk}(t_i)$ with wavelets $\phi_{jk}$. It follows that $y_{jk}$ are independent normal variables with $Ey_{jk} \approx \int f \phi_{jk}$ and $\mathrm{Var}(y_{jk}) = \varepsilon^2 \equiv \sigma^2/N$. See Donoho and Johnstone (1994a, 1995) for details.

Although the variance $\varepsilon^2$ can be fully identified, that is, estimated without error, based on data in (1.1) or (1.6) for square summable $\beta$, that is, $\int_0^1 f^2 < \infty$, implementation of GEB estimators in the nonparametric regression model (4.1) requires an estimation of the variance $\sigma^2$. Among other methods, estimates of $\sigma^2$ can be constructed from observations at the highest resolution level, for example,

$$(4.3) \qquad \hat{\sigma} \equiv \mathrm{MAD}(\sqrt{N}y_{[J]}) \equiv \frac{\mathrm{median}(\sqrt{N}|y_{Jk}|:1 \leq k \leq 2^J)}{\mathrm{median}(|N(0,1)|)},$$

which converges to $\sigma$ at the rate $N^{-\alpha p/(p+1)} + N^{-1/2}$ in Besov balls. The regression function $f$ is then estimated by

$$(4.4) \qquad \hat{f}(i/N) \equiv \sqrt{N} \sum_{j=-1}^{J} \sum_{k=1}^{2^j \vee 1} W_{jk}(i)\, \hat{\beta}_{jk}^{(\varepsilon)}|_{\varepsilon=\hat{\sigma}/\sqrt{N}}$$

via (4.2), where $\hat{\beta}_{jk}^{(\varepsilon)}$ are as in (2.10) and (2.11).

Now we report some simulation results to illustrate the performance of our GEB estimators. Figure 1 plots four examples of regression functions in Donoho and Johnstone (1994a). Normal errors are added to these functions, with signal-to-noise ratio 7, and the resulting response variables $Y_i$, as in (4.1), are plotted against $t_i = i/N$ in Figure 2, with sample size $N = 2048$. Figure 3 reports the GEB estimates (4.4) based on the data in Figure 2, with $j_* = 6$, $\rho_0 = 0.4$, $\eta_n = 0$ and $b_0 = 2$ in (2.10) and (2.11). Figure 4 reports the reconstructions of these regression functions using SURESHRINK in S-plus [Donoho and Johnstone (1995)], also based on the data in Figure 2. The GEB and SURESHRINK estimates look similar in these examples.

4.2. *Random design.* Now consider (4.1) with i.i.d. uniform $t_i$ in $[0, 1]$. We implement GEB methods with Haar basis and provide their optimality properties.

Let $\mathbb{1}_{j,k}(x) \equiv I\{(k-1)/2^j < x \leq k/2^j\}$. The Haar wavelets are $\phi_{j,k} = \sqrt{2^j}(\mathbb{1}_{j+1,2k-1} - \mathbb{1}_{j+1,2k})$, $j \geq 0$, and $\phi_{-1,1} = 1$, and the corresponding wavelet coefficients are

$$(4.5) \quad \beta_{j,k} \equiv \beta_{j,k}(f) \equiv \int_0^1 f\phi_{j,k} = \begin{cases} (\bar{f}_{j+1,2k-1} - \bar{f}_{j+1,2k})/2^{j/2+1}, & j \geq 0, \\ \bar{f}_{0,1}, & j = -1, \end{cases}$$

where $\bar{f}_{j,k} \equiv 2^j \int_0^1 f\mathbb{1}_{j,k}$. Let $N_{j,k} \equiv \sum_i \mathbb{1}_{j,k}(t_i)$ and $\bar{Y}_{j,k} \equiv \sum_i Y_i \mathbb{1}_{j,k}(t_i)/N_{j,k}$. Define

$$(4.6) \quad y_{j,k} \equiv \frac{\delta_{j,k}(\bar{Y}_{j+1,2k-1} - \bar{Y}_{j+1,2k})}{\sqrt{N}(1/N_{j+1,2k-1} + 1/N_{j+1,2k})^{1/2}}, \qquad j \geq 0, \ y_{-1,1} \equiv \bar{Y}_{0,1},$$

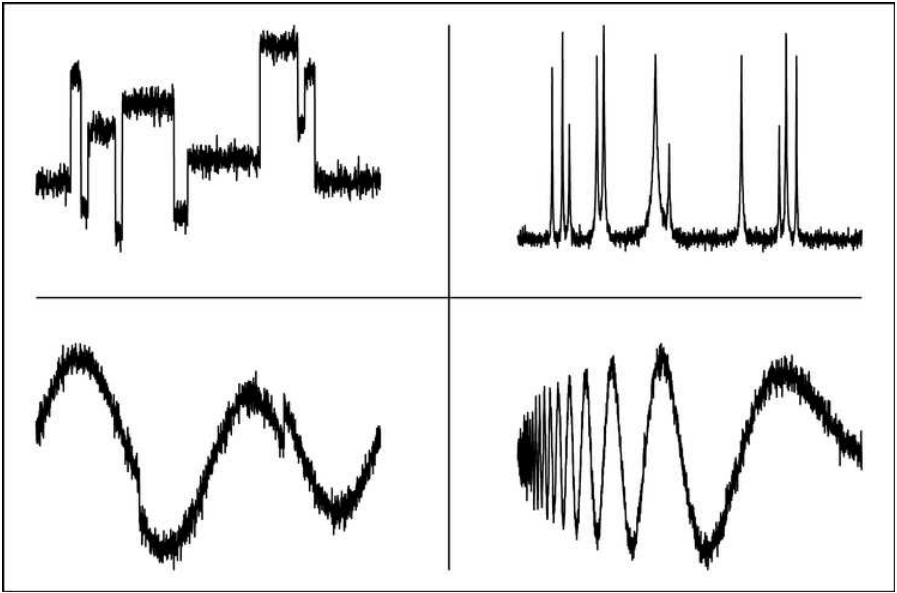

FIG. 2. *Signals + noise with* $N = 2048$; *signal-to-noise ratio is* 7.

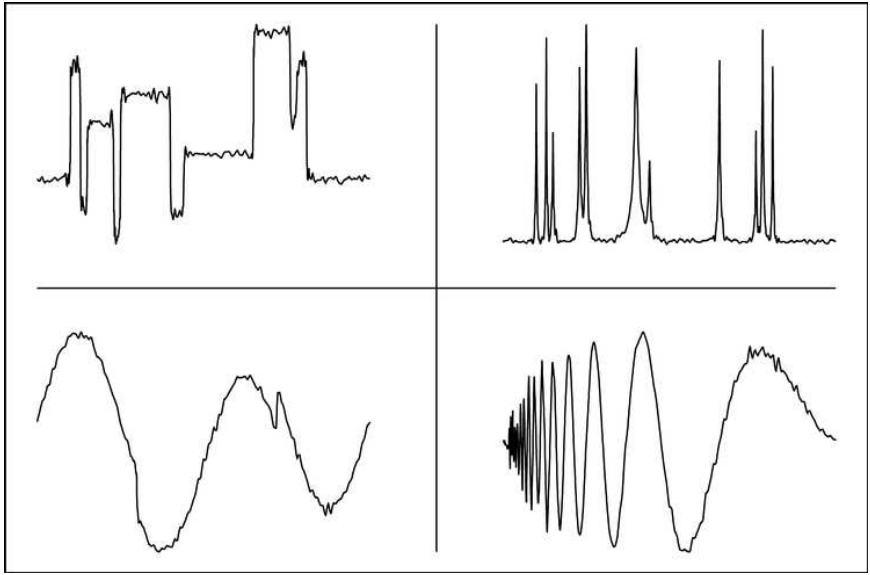

FIG. 3.  *GEB estimate of signals using S8 wavelets; $j_* = 6$, $\rho_0 = 0.4$, $\eta_{[j]} = 0$, $b_0 = 2$.*

where $\delta_{j,k} \equiv I\{N_{j+1,2k-1}N_{j+1,2k} > 0 \text{ or } j = -1\}$. Conditionally on $\{t_i\}$, $y_{j,k}$ are naive estimates of $\beta_{j,k}$ for $\delta_{j,k} = 1$, standardized to have variance $\varepsilon^2 = \sigma^2/N$, and $y_{j,k} \equiv 0$ for $\delta_{j,k} = 0$. In fact, conditionally on $\{t_i\}$, $y_{j,k}$ are inde-

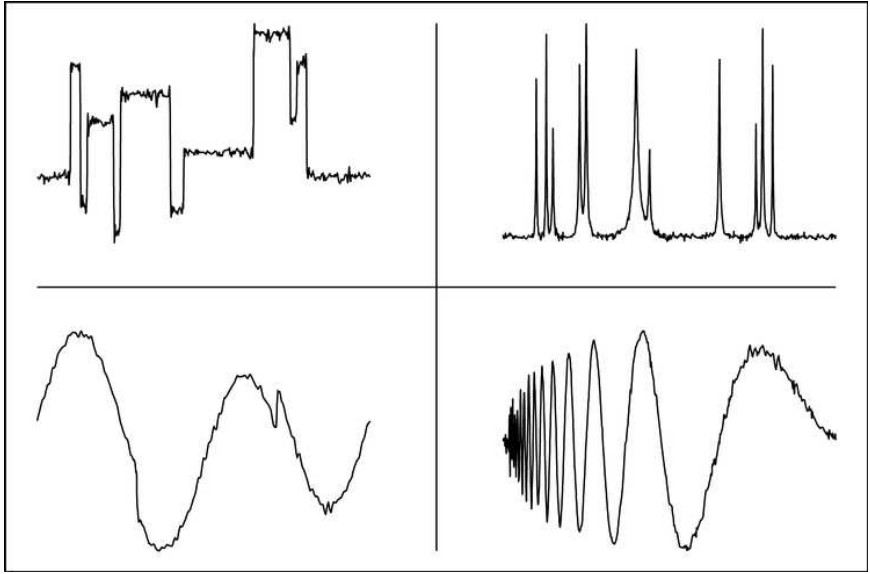

FIG. 4.  *SURESHRINK reconstruction using S8 wavelets.*

pendent $N(\tilde{\beta}_{j,k}, \delta_{j,k}\varepsilon^2)$ variables with

$$
\begin{aligned}
\tilde{\beta}_{j,k} &\equiv \tilde{\beta}_{j,k}(f) \\
&\equiv \frac{\delta_{j,k}(\tilde{f}_{j+1,2k-1} - \tilde{f}_{j+1,2k})}{\sqrt{N}(1/N_{j+1,2k-1} + 1/N_{j+1,2k})^{1/2}}, \qquad j \geq 0, \; \tilde{\beta}_{-1,1} \equiv \tilde{f}_{0,1},
\end{aligned}
\tag{4.7}
$$

where $\tilde{f}_{j,k} \equiv \sum_i f(t_i)\mathbb{1}_{j,k}(t_i)/N_{j,k}$. By the strong law of large numbers, $\tilde{\beta}_{j,k} \to \beta_{j,k}$ as $N \to \infty$.

The statistics $\{y_{j,k}, \delta_{j,k}\}$ in (4.6) are sufficient. Since the data contains no information about $\beta_{j,k}$ for $\delta_{j,k} = 0$, we estimate $\beta_{[j]} \equiv \{\beta_{j,k} : \delta_{j,k} = 1\}$ by GEB based on $y_{[j]} \equiv \{y_{j,k} : \delta_{j,k} = 1\}$,

$$
\hat{\beta}_{j,k} \equiv y_{j,k} I\{j < j_*\} + \delta_{j,k}\varepsilon \hat{t}_{[j]}(y_{j,k}/\varepsilon) I\{j_* \leq j \leq J\},
\tag{4.8}
$$

where $\hat{t}_{[j]}$ is as in (2.8) with $n_j = \sum_k \delta_{j,k}$, $\rho(n)$ and $b(n)$ in (2.11) and

$$
\hat{\varphi}_{[j]}(x) \equiv \frac{\sqrt{2\log n_j}}{n_j} \sum_{k=1}^{2^j} \delta_{j,k} K\left(\frac{x - y_{j,k}/\varepsilon}{(2\log n_j)^{-1/2}}\right),
$$

$$
\hat{\kappa}_{[j]} \equiv 1 - \frac{\sqrt{2}}{n_j} \sum_{k=1}^{2^j} \delta_{j,k} e^{-(y_{j,k}/\varepsilon)^2/2}.
$$

We estimate $f$ by (1.7) via the Parseval identity (1.8).

The following theorem asserts the exactly adaptive minimaxity and super-efficiency of GEB estimators in Besov balls. Let $\bar{f}_j(x) = \sum_{\ell=-1}^{j-1} \sum_k \beta_{\ell,k} \mathbb{1}_{\ell,k}(x)$ be the piecewise average of $f$ at resolution level $j$. For Haar coefficients (4.5), the Besov norm in (3.9) can be written as

$$
\|\beta(f)\|_{p,q}^\alpha = \left\{ |\bar{f}_{0,1}|^q + \sum_{j=0}^\infty 2^{j\alpha q} \left( \int_0^1 |\bar{f}_j - \bar{f}_{j+1}|^p \, dx \right)^{q/p} \right\}^{1/q}.
$$

THEOREM 4.1. *Let $\|f\| \equiv (\int_0^1 f^2)^{1/2}$ and $(\alpha, p)$ satisfy $\alpha^2/(\alpha + 1/2) \geq 1/p - 1/2$. Let $E_f$ be the expectation in (4.1) under which $t_i$ are i.i.d. uniform variables in $(0,1)$. Let $\hat{f} \equiv \hat{f}_N$ be as in (1.7) based on $\hat{\beta} = \{\hat{\beta}_{j,k}\}$ in (4.8), with the cut-off resolution levels $J \equiv J_N$ satisfying $1/\log N \leq \eta_N \equiv 2^{J+1}/N = O(1)$. Then, for all function classes $\mathcal{F} = \{f : \|\beta(f)\|_{p,q}^\alpha \leq C\}$*

$$
\sup_{f \in \mathcal{F}} E_f \|\hat{f}_N - f\|^2 = (1 + \zeta_N) \inf_{\tilde{f}} \sup_{f \in \mathcal{F}} E_f \|\tilde{f} - f\|^2 \sim N^{-2\alpha/(2\alpha+1)}
\tag{4.9}
$$

*with $\zeta_N = o(1)$, provided that $\alpha^2/(\alpha + 1/2) > 1/p - 1/2$ and $\eta_N = o(1)$. Moreover, if $\alpha^2/(\alpha + 1/2) > 1/p - 1/2$ or $\eta_N^{-1} = O(1)$, then (4.9) holds with $\zeta_N = O(1)$ and for all $\|\beta(f)\|_{p,p}^\alpha$-compact classes $\mathcal{F}$*

$$
\sup_{f \in \mathcal{F}} E_f \|\hat{f}_N - f\|^2 = o(1) N^{-2\alpha/(2\alpha+1)}.
\tag{4.10}
$$

For $\delta_{j,k} = 0$, the $N$ observations in (4.1) contain no information about $\beta_{j,k}$ in (4.5). For $2^j \geq N$, this happens for at least half of $\beta_{j,k}$. Thus, the minimax MISE is at least of the order $\max_{f \in \mathcal{F}} \sum_{j,k} \beta_{j,k}^2(f) I\{2^j \geq N\} \sim N^{-2(\alpha+1/2-1/p)}$ in the Besov classes in (4.9). It follows that the condition $\alpha^2/(\alpha+1/2) \geq 1/p - 1/2$, that is, $\alpha + 1/2 - 1/p \geq \alpha/(2\alpha+1)$, is necessary for (4.9). Theorems 4.1 and 9.1 are proved together at the end of the Appendix.

**5. Related problems.** Although the focus of this paper is on the white noise model, our methods have much broader consequences in nonparametric problems and their applications. In addition to the direct implementations in nonparametric regression models in Section 4, the connections between the white noise model and a number of experiments have been recently established in the form of global asymptotic equivalence. This was done by Brown and Low (1996), Donoho and Johnstone (1998) and Brown, Cai, Low and Zhang (2002) for nonparametric regression, by Nussbaum (1996) for the nonparametric density problem and by Grama and Nussbaum (1998) for nonparametric generalized linear models. The impact of such equivalence results is that statistical procedures derived in the white noise model, including those in this paper, can be translated into asymptotically analogous procedures in all other asymptotically equivalent problems. Adaptive estimation in the white noise model (1.1) is also closely related to statistical model selection [cf. Foster and George (1994) and Barron, Birgé and Massart (1999)] and to information theory [cf. Foster, Stine and Wyner (2002)].

There has recently been a spate of papers on adaptive wavelet-based nonparametric methods; see Donoho and Johnstone (1994a, 1995), Donoho, Johnstone, Kerkyacharian and (1995) and Juditsky (1997) on wavelet thresholding in the white noise and nonparametric regression models, Johnstone, Kerkyacharian and Picard (1992) and Donoho, Johnstone, Kerkyacharian and Picard (1996) on related methods in density estimation, Hall, Kerkyacharian and Picard (1998, 1999) and Cai (1999) on block threshold estimators, Abramovich, Benjamini, Donoho and Johnstone (2000) on thresholding based on the false discovery rate, and the recent book of Härdle, Kerkyacharian, Picard and Tsybakov (1998). Adaptive kernel methods were considered by Lepski, Mammen and Spokoiny (1997). These estimators are either nearly adaptive minimax with an extra logarithmic factor in maximum risk in Besov balls (3.9) or rate adaptive for restricted values of $\alpha$ and $p$, for example, $\alpha + 1/2 - 1/p > \{(1/p - 1/2)^+ + \gamma - 1/2\}^+$ in the white noise model, $0 < \gamma < 1/2$, and $\alpha \geq 1/p$ and $p \geq 1$ in nonparametric regression and density problems. This naturally raised the question of the existence of fully rate adaptive estimators for all Besov balls in (3.8), to which Theorem 3.4 provides a positive sharper answer: adaptation to the minimax constants. Cai (2000) pointed out that such sharp adaptation cannot be achieved by separable estimators. The practical value of adaptation for $\alpha < 1/p$ and $p < 1$ is clearly seen from Example 3.1 and Theorem 4.1 and will be further discussed in Section 8. Spatially adaptive methods were also considered by Breiman, Friedman, Olshen and Stone (1984) and Friedman (1991). Johnstone and Silverman (1998, 2004, 2005) proposed a parametric EB approach based on the posterior median for Gaussian errors with respect to a prior as the mixture of the point mass at zero and a given symmetric distribution (e.g., double exponential), with a modified MLE for the mixing probability. Their methods are rate adaptive minimax in all Besov balls and provide stable threshold levels for sparse and dense signals.

Our strategy is to translate high- and infinite-dimensional estimation problems into estimating a sequence of normal means and use block EB methods to derive adaptive estimators. Within each block, one may use general [Robbins (1951, 1956)], linear [Stein (1956), James and Stein (1961) and Efron and Morris (1973)] or other restricted EB methods. From this point of view, the estimator of Efromovich and Pinsker (1984) is block linear EB, while those of Donoho and Johnstone (1995) are block threshold EB. In the wavelet setting, restricted EB could yield exactly adaptive minimax estimators in Besov balls with a fixed primary shape parameter $p$, if $\mathcal{D}_0 \supseteq \{t_{p,c}^* : c > 0\}$, in view of the difference between (3.1) and (3.16), where $t_{p,c}^*$ is the minimax Bayes rule for the class of priors $\{G : \int |\theta|^p \, dG(\theta) \leq c^p\}$. But this is not practical, since the explicit form of $t_{p,c}^*$ is intractable for $p < 2$. In particular, for $p < 2$ the Bayes rules $t_{p,c}^*$ are nonlinear analytic functions, so that linear and threshold estimators do not achieve exact asymptotic minimaxity; see Donoho and Johnstone (1994a, 1998), (3.15) and (7.1) at the resolution level $2^j \sim \varepsilon^{-1/(\alpha+1/2)}$. We further refer to Morris (1983), Robbins (1983) and Berger (1985) for general discussion about EB and Bayes methods.

Adaptive minimax estimation has a number of interpretations. Define

$$\tau(\varepsilon; \hat{\beta}, B) \equiv \frac{\sup_{\beta \in B} R^{(\varepsilon)}(\hat{\beta}, \beta)}{\mathcal{R}^{(\varepsilon)}(B)},$$

where $\mathcal{R}^{(\varepsilon)}(B)$ is the minimax risk in (1.4). Given estimators $\hat{\beta}^{(\varepsilon)}$ and a collection $\mathcal{B}$ of sets $B$ in the parameter space, exactly adaptive minimaxity means $\tau(\varepsilon; \hat{\beta}^{(\varepsilon)}, B) \to 1$ as $\varepsilon \to 0+$ for all $B \in \mathcal{B}$, rate adaptive minimaxity means $\tau(\varepsilon; \hat{\beta}^{(\varepsilon)}, B) = O(1)$, and nearly adaptive minimaxity means that $\tau(\varepsilon; \hat{\beta}^{(\varepsilon)}, B)$ is slowly varying in $\varepsilon$, and with obvious change of notation $\varepsilon \leftrightarrow \sigma/\sqrt{n}$ and $\beta \leftrightarrow f$ for nonparametric regression and density estimation problems. In the wavelet setting, rate and nearly adaptive minimax estimators were derived in Hall and Patil (1995, 1996) and Barron, Birgé and Massart (1999), and block James–Stein estimators were recently investigated by Cavalier and Tsybakov (2001, 2002), in addition to papers cited above. There is a vast literature in nonparametric estimation methods, and asymptotic minimaxity and adaptivity have been commonly used to judge the overall performance of estimators; see comprehensive reviews in Stone (1994), Donoho, Johnstone, Kerkyacharian and Picard (1995) and Barron, Birgé and Massart (1999), and recent books by Efromovich (1999) and Hastie, Tibshirani and Friedman (2001).

**6. Compound estimation of normal means.** Let $(X_k, \theta_k)$, $1 \leq k \leq n$, be random vectors and let $P_{\theta_{(n)}}$ be the conditional probability given $\theta_{(n)} \equiv (\theta_1, \ldots, \theta_n)$ under $P$. Write $P = P_{\theta_{(n)}}$ when $\theta_{(n)}$ is deterministic. Suppose $X_k$ are independent $N(\theta_k, 1)$ variables under the conditional probability $P_{\theta_{(n)}}$. In this section we consider the estimation of $\theta_k$ under the compound squared error loss $n^{-1} \sum_{k=1}^{n} (\hat{\theta}_k - \theta_k)^2$, that is, the estimation of normal means within a single block or resolution level based on (2.2) or (1.1), scaled to the unit variance.

Let $X \sim N(\theta, 1)$ under $P_\theta$. Define the Bayes risks for Borel $t(\cdot)$ and their minimum by

$$(6.1) \quad R(t, G) \equiv \int E_\theta (t(X) - \theta)^2 \, dG(\theta), \qquad R^*(G) \equiv \inf_t R(t, G).$$

As pointed out by Robbins (1951), the compound mean squared error for the use of $\hat{\theta}_k = t(X_k)$ is $R(t, G_n)$, where $G_n$ is the mixture of the marginal distributions of $\theta_{(n)}$,

$$(6.2) \quad G_n(x) \equiv \frac{1}{n} \sum_{k=1}^{n} P\{\theta_k \leq x\}.$$

We drive GEB estimators whose compound risk approximates the ideal Bayes risk $R^*(G_n)$. We measure the performance of this ideal approximation via oracle inequalities of the form

$$(6.3) \quad \frac{1}{n} \sum_{k=1}^{n} E(\hat{\theta}_k - \theta_k)^2 - R^*(G_n) \leq r(n, G_n),$$

where $r(n, G)$ are functionals of $n$ and univariate distributions $G$ only. The definition of $r(n, G)$ may vary in different statements in the sequel, as long as (6.3) holds under specified conditions.

The components of the vector $\theta_{(n)}$ are assumed to be independent in Theorem 6.1 below. In all other theorems, conditions on $\theta_{(n)}$ are imposed only through the mixture $G_n$ in (6.2), so that $\theta_k$ are allowed to be stochastically dependent. The independence assumption on $\theta_{(n)}$ in Theorem 6.1 accommodates the two important special cases of deterministic and i.i.d. $\{\theta_k\}$. This allows us to apply Theorem 6.1 conditionally on $\theta_{(n)}$ whenever $r(n, G)$ in (6.3) is concave in $G$. Note that if $r(n, G)$ is concave in $G$, (6.3) follows from its conditional version given $\theta_{(n)}$, since $R^*(G)$ is always concave in $G$ due to the linearity of $R(t, G)$ in $G$ in (6.1).

6.1. *GEB estimators.* Zhang (1997) proposed the following GEB estimators:

(6.4) $\quad \hat{\theta}_k \equiv \hat{t}_{n,\rho}(X_k), \qquad k \leq n, \qquad \hat{t}_{n,\rho}(x) \equiv x + \dfrac{\hat{\varphi}'_n(x)}{\max(\hat{\varphi}_n(x), \rho)},$

where $\rho \equiv \rho_n \to 0+$, $1/n \leq \rho < 1/\sqrt{2\pi}$, and $\hat{\varphi}_n$ is the kernel estimator

(6.5) $\quad \hat{\varphi}_n(x) \equiv \dfrac{1}{n}\sum_{k=1}^{n} a_n K(a_n(x - X_k)) = \int_{-a_n}^{a_n} \dfrac{e^{-ixu}}{2\pi} \sum_{k=1}^{n}\dfrac{e^{iuX_k}}{n}\, du$

with the kernel $K(x) \equiv \sin(x)/(\pi x)$. We use the special $a_n = \sqrt{2\log n}$ throughout the sequel, which provides the best bounds in this paper. We first describe an improved version of the oracle inequality of Zhang (1997) and its immediate consequences.

THEOREM 6.1. *Suppose the components of $\theta_{(n)} \equiv (\theta_1, \ldots, \theta_n)$ are independent variables. Let $\hat{\theta}_k \equiv \hat{t}_{n,\rho}(X_k)$ be the GEB estimator in (6.4) with $\rho^{-1}(\log n)^{1/4}/\sqrt{n} = o(1)$. Then (6.3) holds with*

(6.6) $\quad r(n, G) = \Delta(\rho, G) + \{1 + \eta(n, \rho)\}\Delta^*(n, \rho),$

*where $\eta(n, \rho) = o(1)$ depending on $(n, \rho)$ only,*

(6.7) $\quad \Delta(\rho, G) \equiv \int_{-\infty}^{\infty} \{\varphi'_G/\varphi_G\}^2 \{1 - \varphi_G/(\varphi_G \vee \rho)\}^2 \varphi_G\, dx$

*with $\varphi_G \equiv \varphi(x; G)$ in (2.5), $G_n$ in (6.3) is as in (6.2), and*

(6.8) $\quad \Delta^*(n, \rho) \equiv \{\sqrt{(2/3)\log n} + \sqrt{-\log(\rho^2)}\}^2 \dfrac{\sqrt{2\log n}}{\pi \rho n}.$

REMARK. (i) The oracle inequality (6.6) was proved in Zhang [(1997), Theorem 1] under the stronger condition $\rho^{-1}\sqrt{(\log n)/n} = o(1)$. The weaker condition is needed since (2.11) is used in this paper. (ii) By (6.7), $R^*(G) + \Delta(\rho, G) \leq 1$, since

(6.9) $\quad R^*(G) = 1 - \int \left(\dfrac{\varphi'_G}{\varphi_G}\right)^2 \varphi_G\, dx.$

The main consequences of Theorem 1 of Zhang (1997) and Theorem 6.1 above under weaker conditions on $\rho_n$ are *asymptotic minimaxity* and *asymptotic optimality*. It is well known that the minimax mean squared error for compound estimation of normal means is the common variance.

THEOREM 6.2. *Let $\hat{\theta}_k \equiv \hat{t}_{n,\rho}(X_k)$ be as in (6.4) with $\rho \equiv \rho_n \to 0$ and $(\log n)^{1/4}/(\rho\sqrt{n}) \to 0$.*

(i) *Asymptotic minimaxity:* For the $\Delta^*(n,\rho)$ in (6.8),

$$\sup_{\theta_{(n)}} \frac{1}{n}\sum_{k=1}^{n} E_{\theta_{(n)}}(\hat{\theta}_k - \theta_k)^2 \leq 1 + (1+o(1))\Delta^*(n,\rho_n) \to 1.$$

(ii) *Asymptotic optimality:* If $G_n$ converges in distribution, then

(6.10) $$\frac{1}{n}\sum_{k=1}^{n} E(\hat{\theta}_k - \theta_k)^2 - R^*(G_n) \to 0,$$

where $R^*(G)$ and $G_n$ are as in (6.1) and (6.2). Moreover, for $m_n = o(1/\rho_n)$ and any stochastically bounded family $\mathcal{G}$ of distributions, (6.10) holds if for certain $0 \leq w_{n,0} \to 0$ and distributions $H_{n,0}$, $G_n(x) = \sum_{j=0}^{m_n} w_{n,j} H_{n,j}(x - c_{n,j})$ with $H_{n,j} \in \mathcal{G}$ for $j \geq 1$ and reals $w_{n,j} \geq 0$ and $c_{n,j}$, that is, $G_n$ are within $o(1)$ mass from mixtures of at most $o(1/\rho_n)$ arbitrary translations of distributions in $\mathcal{G}$. In particular, (6.10) holds if $\int_{|x-c_n|>m_n} dG_n(x) \to 0$ for $m_n = o(1/\rho_n)$ and certain constants $c_n$.

REMARK. (i) Zhang [(1997), Proposition 2 and Corollary 3] pointed out that (6.10) holds when $G_n$ converges in distribution or when $G_n$ are arbitrary discrete distributions with no more than $o(1/\rho_n)$ components. The weaker condition in Theorem 6.2(ii) is equivalent to $G_n(A_n) \to 1$ for certain unions $A_n$ of at most $m_n = o(1/\rho_n)$ intervals of unit length. This demonstrates the extent of adaptivity of GEB estimators when $\{\theta_k\}$ has many clusters.

(ii) The proof of Theorem 6.2(ii) utilizes the following inequality: for all distributions $H_j$ and weights $w_j > 0$ with $\sum_{j=0}^{m} w_j = 1$,

(6.11) $$G = \sum_{j=0}^{m} w_j H_j \implies \Delta(\rho, G) \leq \sum_{j=0}^{m} w_j \Delta(\rho/w_j, H_j).$$

(iii) The locally uniform asymptotic optimality criterion in (6.10) is slightly stronger than the usual one for fixed $G = G_n$ in the EB setting.

6.2. *Oracle inequalities based on tail-probabilities and moments.* We shall derive more explicit oracle inequalities in terms of the tail and moments of $G_n$. Define

(6.12) $$\overline{G}(x) \equiv \int_{|u|>x} dG(u),$$

$$\mu_p(G) \equiv \left(\int |u|^p \, dG(u)\right)^{1/p}, \qquad 0 < p < \infty.$$

LEMMA 6.1. *Let $x \geq 0$, $0 < \rho < 1/\sqrt{2\pi}$ and $\varphi_G$ be as in (2.5). Then*

$$\Delta(\rho, G) \leq \int_{\varphi_G \leq \rho} \left(\frac{\varphi'_G}{\varphi_G}\right) \varphi_G$$

(6.13)

$$\leq \overline{G}(x) + 2x\rho \max\{\widetilde{L}^2(\rho), 2\} + 2\rho\sqrt{\widetilde{L}^2(\rho) + 2},$$

*where $\Delta(\rho, G)$ is as in (6.7) and $\widetilde{L}(\rho) \equiv \sqrt{-\log(2\pi\rho^2)}$. Furthermore, for $x = \widetilde{L}(\rho)/2$,*

(6.14) $\qquad \Delta(\rho, G) \leq \overline{G}(x) + \overline{G}^2(x)(1 - \overline{G}(x)) + 2\rho\sqrt{\widetilde{L}^2(\rho) + 2}.$

Lemma 6.1 is used in combination with Theorem 6.1 to produce more explicit oracle inequalities in Theorems 6.3 and 6.4 below, with (6.13) for stochastically large $G_n$ and (6.14) for stochastically small $G_n$. For stochastically very small $G_n$ and $-\log \rho_n^2 \leq (1 + o(1)) \log n$, the leading term in the combination of (6.6) and (6.14) is

(6.15) $\quad \Delta^*(n, \rho_n) + 2\rho_n\sqrt{\widetilde{L}(\rho_n) + 2} \leq (1.724 + o(1))\frac{\log n}{\sqrt{n}}\left(\frac{\rho_n^*}{\rho_n} + \frac{\rho_n}{\rho_n^*}\right)$

with equality for $-\log \rho_n^2 = (1 + o(1)) \log n$, where $\rho_n^* \equiv 0.6094\sqrt{2(\log n)/n}$. The choice of $\rho \approx \rho_n^*$ and the oracle inequalities below are not necessarily optimal, since crude bounds are used at several places in the proofs. In principle, we may use data-driven $\rho$ via any methods of choosing tuning parameters, but this is beyond the scope of this paper. In what follows, we denote by $\eta_n$ constants depending on $n$ only and satisfying $\eta_n \to 0$.

THEOREM 6.3. *Let $\hat{\theta}_k$ be as in (6.4) with $\rho \equiv \rho(n)$ in (2.11). Then (6.3) holds with*

(6.16)
$$r(n, G) = \inf_x \left\{\overline{G}(x) + (1 + \eta_n)x\sqrt{8}\rho_0\frac{(\log n)^{3/2}}{\sqrt{n}}\right\} + (1 + \eta_n)C_0(\rho_0)\frac{\log n}{\sqrt{n}}$$

$$\leq (1 + \eta_n)\left[C_p\left(\frac{\rho_0\mu_p(G)(\log n)^{3/2}}{\sqrt{n/8}}\right)^{p/(p+1)} + C_0(\rho_0)\frac{\log n}{\sqrt{n}}\right],$$

*where $\overline{G}(x)$ and $\mu_p(G)$ are as in (6.12), $C_p \equiv p^{1/(p+1)} + p^{-p/(p+1)} \leq 2$ and $C_0(x) \equiv 1.724(0.6094/x + x/0.6094)$. Moreover, for $x_n \equiv \sqrt{-\log(2\pi\rho_n^2)}/2$, inequality (6.3) holds with*

(6.17) $\qquad r(n, G) = (5/4)\overline{G}(x_n) + (1 + \eta_n)C_0(\rho_0)(\log n)/\sqrt{n}.$

Theorem 6.3 provides the asymptotic optimality of GEB estimators with convergence rates $\{(\log n)^{3/2}/\sqrt{n}\}^{-p/(p+1)}$ in (6.10) for dependent $\{\theta_k\}$ with bounded $\mu_p(G_n)$.

6.3. *Stochastically very small distributions and threshold estimators.* The risk bounds in Theorems 6.1 and 6.3 are not very useful if an overwhelming majority of $\theta_k$ are essentially zero, for example, $\mu_2(G_n) \leq 1/\sqrt{n}$. For these stochastically very small empirical distributions $G_n$, threshold estimators may outperform the GEB estimators (6.4).

Soft threshold estimators are defined by

$$\hat{\theta}_k \equiv s_\lambda(X_k), \qquad s_\lambda(x) \equiv \text{sgn}(x)(|x| - \lambda)^+, \tag{6.18}$$

where $\lambda > 0$ is a threshold level. Hard threshold estimators are defined by functions $h_\lambda(x) \equiv xI\{|x| > \lambda\}$. Hard and soft threshold estimators have similar properties. We consider soft threshold estimators so that sharp oracle inequalities in Lemma 6.2 below can be utilized.

The performance of (6.18) is commonly compared with $\kappa(G_n) = E\sum_{k=1}^n (\theta_k^2 \wedge 1)/n$ given in (2.7). For $A \subseteq \{1, \ldots, n\}$, let $t_A$ be the estimator defined by $\hat{\theta}_k = X_k I\{k \in A\}$. Since the MSE of $X_k$ is smaller than the MSE of $\hat{\theta}_k = 0$ iff $|\theta_k| > 1$, $\kappa(G_n) = \inf_A R(t_A, G_n)$ when $\theta_{(n)}$ is deterministic. Thus, $\kappa(G_n)$ is the ideal risk for a different oracle expert, someone with the knowledge of the best choice of $A$, who always uses the best $t_A$.

LEMMA 6.2. *Let $s_\lambda \equiv s_\lambda(x)$ be as in (6.18) and let the risk $R(t, G)$ be as in (6.1). Then $E_\theta(s_\lambda(X) - \theta)^2$ is increasing in $|\theta|$, and*

$$R(s_\lambda, G) \leq \int \min\left\{u^2 + \frac{4}{\lambda^3}\varphi(\lambda), \lambda^2 + 1\right\} dG(u). \tag{6.19}$$

*Consequently, for $\lambda = \sqrt{2\log n}$ and with $\mu_p(G)$, $\overline{G}(x)$ and $\kappa(G)$ as in (6.12) and (2.7), $p \leq 2$,*

$$R(s_\lambda, G) \leq \min\left\{\frac{\mu_p^p(G)}{(2\log n + 1)^{p/2-1}}, (2\log n)\overline{G}(1) + \kappa(G)\right\} + \frac{\sqrt{2}}{n(\log n)^{3/2}}. \tag{6.20}$$

The inequalities in Lemma 6.2 are essentially the oracle inequality of Donoho and Johnstone (1994a). The improvement with the extra factor $(\log n)^{-3/2}$ in the second term on the right-hand side of (6.20) is needed when we apply it to all high-resolution levels $j$ near the infinity in the sequence model (1.1). Lemma 6.2 implies $R(s_\lambda, G) \leq (\lambda^2 + 1)\kappa(G) + 4\lambda^{-3}\varphi(\lambda)$, which is an oracle inequality since it compares $R(s_\lambda, G)$ with the ideal risk $\kappa(G)$. For $\lambda = \sqrt{2\log n}$, Foster and George (1994) showed that $\lambda^2 + 1$ is the optimal *risk inflation factor* from a model selection point of view.

Since $\overline{G}(x) \leq \kappa(G) \leq \mu_p^p(G)$ for $p \leq 2$ and $x \geq 1$, the GEB oracle inequality (6.17) (with $x_n \to \infty$) can be directly compared with (6.20). The risk bound for the threshold estimator is of larger order than the regret of the GEB estimator if $\kappa(G_n)\sqrt{n}/\log n \to \infty$.

6.4. *Hybrid GEB methods.* In the white noise model (1.1), $\sum_{j,k} \beta_{jk}^2 < \infty$, so that the ideal risk $\kappa(G_n)$ converges to zero as $n = 2^j \to \infty$. Thus, the performance of the GEB estimator (6.4) could be enhanced if hybrid estimators are used, that is, switching to the threshold estimator (6.18) for small $\kappa(G_n)$. By Zhang (1990), $\overline{G}_n(x)$ and thus $\kappa(G_n) = \int_0^1 \overline{G}_n(u) \, du^2$ can be estimated only at logarithmic rates. Our strategy is to construct hybrid estimators based on accurate estimates of the order of $\kappa(G_n)$.

The order of magnitude of $\kappa(G)$ in (2.7) is the same as that of

$$(6.21) \quad \tilde{\kappa}(G) \equiv 1 - \int \sqrt{2} e^{-x^2/2} \varphi(x; G) \, dx = 1 - \int \exp(-u^2/4) \, dG(u).$$

In fact, since $(1 - 1/e)x \le 1 - e^{-x} \le x$ for $0 \le x \le 1$,

$$(6.22) \quad \begin{aligned} \frac{e-1}{4e} \kappa(G) &\le (1 - 1/e) \int \left(\frac{u^2}{4} \wedge 1\right) dG(u) \\ &\le \tilde{\kappa}(G) \le \int \left(\frac{u^2}{4} \wedge 1\right) dG(u) \le \kappa(G). \end{aligned}$$

Thus, the order of $\kappa(G_n)$ can be estimated by

$$(6.23) \quad \hat{\kappa}_n \equiv 1 - \frac{1}{n} \sum_{k=1}^n \sqrt{2} \exp(-X_k^2/2).$$

This suggests the following hybrid estimators:

$$(6.24) \quad \hat{\theta}_k \equiv \tilde{t}_n(X_k), \qquad \tilde{t}_n(x) \equiv \hat{t}_{n,\rho,\lambda,b}(x) \equiv \begin{cases} \hat{t}_{n,\rho}(x), & \text{if } \hat{\kappa}_n > b, \\ s_\lambda(x), & \text{if } \hat{\kappa}_n \le b, \end{cases}$$

where $\hat{t}_{n,\rho}(\cdot)$, $s_\lambda(\cdot)$ and $\hat{\kappa}_n$ are as in (6.4), (6.18) and (6.23), respectively. For definiteness, we choose in (6.24) $\rho \equiv \rho(n)$ and $b \equiv b(n)$ in (2.11) and $\lambda = \sqrt{2 \log n}$, unless otherwise stated. This choice of $\rho_n$ optimizes the order of risk bound (6.17). The choice of $\lambda_n$ matches the universal thresholding [Donoho and Johnstone (1994a)] and provides the optimal risk inflation factor [Foster and George (1994)]. The choice $b_n$ here ensures the use of (6.4) for large $\kappa(G_n)\sqrt{n}/\log n$.

LEMMA 6.3. *Suppose that $\{\theta_k\}$ are independent variables under the expectation $E$. Let $\hat{t}_n \equiv \hat{t}_{n,\rho,\lambda,b}$ be the hybrid estimator in (6.24) with $\lambda \equiv \lambda_n =$*

$\sqrt{2\log n}$. Then

$$\sum_{k=1}^{n} \frac{E(\hat{t}_n(X_k) - \theta_k)^2}{n}$$

(6.25)
$$\leq \begin{cases} \sum_{k=1}^{n} E(\hat{t}_{n,\rho}(X_k) - \theta_k)^2/n + (2 + \eta_n)(\log n)/n^2, & \tilde{\kappa}(G_n) \geq b_n^+, \\ R(s_\lambda, G_n) + (1 + \eta_n)(\log n)^2/(\pi^2 \rho^2 n^3), & \tilde{\kappa}(G_n) \leq b_n^-, \\ \sum_{k=1}^{n} E(\hat{t}_{n,\rho}(X_k) - \theta_k)^2/n + R(s_\lambda, G_n), & \text{otherwise,} \end{cases}$$

with $\eta_n \to 0$ uniformly for all choices of $\rho \equiv \rho_n$ and $b \equiv b_n$, where $\hat{t}_{n,\rho}$, $s_\lambda$ and $\tilde{\kappa}(G)$ are as in (6.4), (6.18) and (6.21), respectively, $b_n^+ \equiv b_n + \sqrt{2(\log n)/n}$, and $b_n^- \equiv b_n - \sqrt{3(\log n)/n}$.

REMARK. Let $(\rho, b)$ be as in (2.11) and $\lambda = \sqrt{2\log n}$. By (6.17) and the fact that $\overline{G}(1) \leq \kappa(G)$, (6.3) holds with $r(n, G) = O(1)(\log n)/\sqrt{n}$ for the GEB estimator (6.4) when $b_n^- < \tilde{\kappa}(G_n) < b_n^+$.

Theorem 6.4 below provides oracle inequalities for (6.24) in terms of the tail of $G_n$ in (6.2).

THEOREM 6.4. *Let $\hat{\theta}_k = \hat{t}_{n,\rho,\lambda,b}(X_k)$ be the hybrid GEB estimator in (6.24) with $(\rho, b)$ in (2.11) and $\lambda = \sqrt{2\log n}$. Then there exists a constant $M < \infty$ such that (6.3) holds with*

(6.26)
$$r(n, G) = M \min\{r_0(n, G), r_{p \wedge 2}(n, \mu_p(G))\} + \frac{1 + \eta_n}{n(\log n + 1)^{3/2}} \quad \forall p > 0, n,$$

*where $r_p(n, C)$ and $\mu_p(G)$ are as in (3.3) and (6.12), and with $\overline{G}$ as in (6.12),*

(6.27)
$$r_0(n, G) \equiv \min\left(1, \int_0^{\log n} \overline{G}(\sqrt{u})\, du, \right. \\ \left. \frac{(\log n)^2}{\sqrt{n}} + \inf_{x \geq 1}\left[\overline{G}(x) + x\frac{(\log n)^{3/2}}{\sqrt{n}}\right]\right).$$

REMARK. It follows from our proof (with slight modification) that if larger $\lambda = \sqrt{2(1 + A_0)\log n}$ is used in (6.24) with $A_0 > 0$, Theorem 6.4 holds with $(1 + \eta_n)/\{n^{1+A_0}(\log n)^{3/2}\}$ as the second term on the right-hand side of (6.26). See the remark below Theorem 3.2.

Theorem 6.4 implies that the compound risk is approximately $n^{-1}(\log n)^{-3/2}$ when $\theta_k = 0$ for all $k$. Proposition 6.1 facilitates applications of Theorem 6.4.

PROPOSITION 6.1. *Let $r_0(n,G)$ and $r_p(n,C)$ be as in (6.27) and (3.3). Let $w' \wedge w'' > 0$.*

(i) *$r_0(n,G)$ is concave in $G$ and $r_0(n,G) \leq 3 r_{p \wedge 2}(n, \mu_p(G))$ for all $p > 0$.*
(ii) *If $\overline{G} \leq w'\overline{G}' + w''\overline{G}''$ for two distributions $G'$ and $G''$, then $r_0(n,G) \leq r_0(n; w'G') + r_0(n; w''G'')$.*

6.5. *Minimax risks in $\ell_p$ balls.* Now we compare the minimax risk

$$(6.28) \quad \mathcal{R}_n(\Theta) \equiv \inf_{\hat{\theta}_{(n)}} \sup_{\theta_{(n)} \in \Theta} \frac{1}{n} \sum_{k=1}^n E_{\theta_{(n)}} (\hat{\theta}_k - \theta_k)^2, \quad \Theta \subset \mathbb{R}^n,$$

in $\ell_p$ balls with the maximum of the Bayes risk $R^*(G)$ (6.1) in $L_p$ balls. Here $\theta_{(n)} \equiv (\theta_1, \ldots, \theta_n)$ are considered as deterministic vectors and the minimization is taken over all estimators $\hat{\theta}_{(n)}$ based on $X_{(n)}$. Our result is based on Proposition 6.2 below, which provides the continuity of the Bayes risk $R^*(G)$ in $G$. Let $\|x_{(n)}\|_{p,n} \equiv (\sum_{k=1}^n |x_k|^p)^{1/p}$ as in (3.9). The $\ell_p$ balls are defined as

$$\Theta_{p,n}(C) \equiv \{\theta_{(n)} : n^{-1/p} \|\theta_{(n)}\|_{p,n} \leq C\},$$

while the $L_p$ balls are $\{G : \mu_p(G) \leq C\}$ with the $\mu_p(G)$ in (6.12).

PROPOSITION 6.2. *Let $R^*(G)$ be as in (6.1). For all distributions $H_1$ and $H_2$ in $\mathbb{R}$*

$$(6.29) \quad |R^*((1-w)H_1 + wH_2) - R^*(H_1)| \leq w\{1 + \sqrt{2\log(\sqrt{2}/w)}\}^2.$$

*Furthermore, if there exist random variables $\tilde{\theta}_k \sim G_k$ with $P\{|\tilde{\theta}_1 - \tilde{\theta}_0| \leq \eta_2\} \geq 1 - \eta_1 \geq 0$, then*

$$(6.30) \quad |R^*(G_1) - R^*(G_0)| \leq 2\eta_1 \{1 + \sqrt{2\log(\sqrt{2}/\eta_1)}\}^2 + \sqrt{8}\left\{1 + \frac{1}{\sqrt{\pi}}\right\}\eta_2.$$

PROPOSITION 6.3. *Let $p' \equiv p \wedge 2$ and $\Psi(u) \equiv \{u \log_+^2 (1/u)\}^{1/3}$. Define*

$$(6.31) \quad r_p^*(n,C) \equiv \begin{cases} \Psi(\log_+(C)/(np^2)), & \text{if } C \geq 1, \\ \Psi(C^{2p'}\{\log_+(1/C^{p'})\}^{2-p'}/(np^2 p')), & \text{if } C < 1, \end{cases}$$

*and $r_\infty^*(n,C) \equiv 0$. Then there exists a universal constant $M$ such that for all $0 < p \leq \infty$*

$$(6.32) \quad \mathcal{R}_n(\Theta_{p,n}(C)) \leq \sup_{\mu_p(G) \leq C} R^*(G) \leq \mathcal{R}_n(\Theta_{p,n}(C)) + M r_p^*(n,C).$$

REMARK. Let $\lambda \equiv \{2\log_+(1/C^{p'})\}^{1/2}$. By (6.19) of Lemma 6.2, uniformly in $p$ as $C^{p'} \to 0$,

$$\sup_{\mu_p(G) \leq C} R^*(G) \leq C^{p'}(\lambda^{2-p'} + 1 + 4/\lambda^3)$$
$$(6.33) \qquad\qquad = (1 + o(1))C^{p'}\{2\log_+(1/C^{p'})\}^{1-p'/2}.$$

Thus, for small $C^{p'}$, (6.32) is sharp only when $Mr_p^*(n, C)$ is smaller than the right-hand side of (6.33), that is, large $C^{p'}(np^2p')/\{\log(np^2p')\}^{1+p'/2}$.

The minimax risk in $\ell_p$ balls and the maximum Bayes risk in $L_p$ balls have been studied by Donoho and Johnstone (1994b), who proved

$$\lim_{C \to 0+} \frac{\sup_{\mu_p(G) \leq C} R^*(G)}{C^{p \wedge 2}\{-2\log(C^{p \wedge 2})\}^{(1-(p \wedge 2)/2)}} = 1 \qquad \forall p > 0,$$

$$(6.34) \qquad \mathcal{R}_n(\Theta_{p,n}(C)) \leq \sup_{\mu_p(G) \leq C} R^*(G)$$
$$\leq b^2 \sup_{\mu_p(G) \leq C/b} R^*(G) \qquad \forall p > 0, \ b \geq 1,$$

and under the extra condition $C^p n/(\log n)^{p/2} \to \infty$ for $p < 2$, $\mathcal{R}_n(\Theta_{p,n}(C)) \approx \sup_{\mu_p(G) \leq C} R^*(G)$ as $n \to \infty$. Proposition 6.3 is derived from Proposition 6.2, (6.34) and Lemma A.3 in the Appendix.

6.6. *Adaptive minimax estimation in $\ell_p$ balls.* An immediate consequence of Theorem 6.4 is the adaptive minimaxity of the GEB estimators in $\ell_p$ balls $\Theta_{p,n}(C)$, in view of the result of Donoho and Johnstone (1994b) on the equivalence of the minimax risk in $\ell_p$ balls and the maximum Bayes risk in $L_p$ balls.

THEOREM 6.5. *Let $\mathcal{R}_n(\Theta)$ be the minimax risk in (6.28) and $\hat{\theta}_{(n)} \equiv (\hat{\theta}_1, \ldots, \hat{\theta}_n)$ be the hybrid GEB estimator in Theorem 6.4. If $C^p\sqrt{n}/(\log n)^{1+(p \wedge 2)/2} \to \infty$, then*

$$(6.35) \qquad \sup_{\theta_{(n)} \in \Theta_{p,n}(C)} \frac{1}{n}\sum_{k=1}^n E(\hat{\theta}_k - \theta_k)^2 = \{1 + o(1)\}\mathcal{R}_n(\Theta_{p,n}(C)).$$

*Moreover, if $C^p n/(\log n)^{(p \wedge 2/2)} \to \infty$, then (6.35) holds with the $o(1)$ replaced by $O(1)$.*

**7. Oracle and risk inequalities for block GEB estimators.** We provide here stronger versions of the theorems in Section 3. This is accomplished by inserting the inequalities in Section 6 in individual blocks or resolution levels. Throughout the section, $E_\beta^{(\varepsilon)}$ denotes the expectation of models (2.2) or (1.1) and $\beta$ is treated as a deterministic sequence. Performance of GEB estimators in more general classes of $\beta$ will be considered in Section 8.

7.1. *Oracle inequalities.* Consider the general sequence (2.2). It follows from (3.1) that

$$(7.1) \quad R^{(\varepsilon,*)}(\beta) = \varepsilon^2 \sum_j n_j \min_{t(\cdot)} R(t, G_{[j]}^{(\varepsilon)}) = \varepsilon^2 \sum_j n_j R^*(t, G_{[j]}^{(\varepsilon)}),$$

where $R(t,G)$ and $R^*(G)$ are as in (6.1) and $G_{[j]}^{(\varepsilon)}$ are as in (2.4). By (2.3) and (2.8)

$$(7.2) \quad \begin{aligned} R^{(\varepsilon)}(\hat{\beta}^{(\varepsilon)}, \beta) &= \sum_j \sum_{k \in [j]} E_\beta^{(\varepsilon)} (\hat{\beta}_k^{(\varepsilon)} - \beta_k)^2 \\ &= \sum_j \varepsilon^2 \sum_{k \in [j]} E_\beta^{(\varepsilon)} (\hat{t}_{[j]}(X_k) - \theta_k)^2, \end{aligned}$$

where $(X_k, \theta_k) \equiv (y_k, \beta_k)/\varepsilon$. Since (2.8) is the implementation of (6.24) in block $j$, application of Theorem 6.4 in individual blocks in (7.2) and (7.1) yields Theorem 3.1 and the following theorem.

THEOREM 7.1. *Let $\hat{\beta}^{(\varepsilon)}$ and $G_{[j]}^{(\varepsilon)}$ be as in Theorem 3.1, let $R^{(\varepsilon,*)}(\beta)$ be as in (3.1) and let $r_0(n,G)$ be as in (6.27). Then there is a universal constant $M < \infty$ such that*

$$(7.3) \quad R^{(\varepsilon)}(\hat{\beta}^{(\varepsilon)}, \beta) - R^{(\varepsilon,*)}(\beta) \leq M\varepsilon^2 \sum_j \left\{ n_j r_0(n_j, G_{[j]}^{(\varepsilon)}) + \frac{1}{(\log n_j + 1)^{3/2}} \right\}.$$

7.2. *Uniform ideal adaptation in Besov balls.* In the wavelet setting (1.1), $n_j = 2^j$, and $\|\beta\|_{p,q}^\alpha \leq C$ iff for certain $C_j \geq 0$ with $(\sum_j C_j^q)^{1/q} = C$,

$$(7.4) \quad \mu_p(G_{[j]}^{(\varepsilon)}) = \left( \frac{1}{2^j} \sum_{k=1}^{2^j \vee 1} \left|\frac{\beta_{jk}}{\varepsilon}\right|^p \right)^{1/p} \leq 2^{-j(\alpha+1/2)} \frac{C_j}{\varepsilon} \leq 2^{-j(\alpha+1/2)} \frac{C}{\varepsilon} \quad \forall j,$$

in view of (6.12), (2.4) and (3.9). Thus, the bound in Theorem 3.1 can be explicitly computed to provide uniform convergence rates for the regret of the GEB estimator (2.10).

We first define certain constants and bounded nonincreasing slowly varying functions. Set

$$(7.5) \quad \alpha_1 \equiv 2\alpha + 1/2 - 1/p', \qquad \alpha_2 \equiv \min(\alpha_1, \alpha + 1/2),$$

with $p' \equiv \min(p, 2)$, and

$$
\gamma_1 \equiv 2 - \frac{1/2 + 1/p'}{\alpha_1 + 1/2},
$$
(7.6)
$$
\gamma_2 \equiv \begin{cases} 1 + 3/(2\alpha + 2), & \text{if } \alpha p' > 1, \\ 3 - (1/2 + 2/p')/(\alpha_1 + 1/2), & \text{otherwise.} \end{cases}
$$

Let $\gamma \equiv p'(\alpha + 1/2) - 1$ and $\tilde{p} \equiv 1/(1 - p'/q)^+ \in [1, \infty]$. Define

$$
(7.7) \quad L^{(1)}(\varepsilon) \equiv L^{(1)}_{\alpha,p,q}(\varepsilon) \equiv (1 + \alpha)^{-\gamma_1} \left[ c'_{\alpha,p,q} + \left( \frac{\log_+(1/\varepsilon)}{\alpha + 1} \right)^{p'/2 - 1} c''_{\alpha,p,q} \right]
$$

with $c'_{\alpha,p,q} \equiv (1 - 1/2^{\tilde{p}\gamma})^{-1/\tilde{p}}$ and $c''_{\alpha,p,q} \equiv (1 - 1/2^{\tilde{p}\gamma})^{-1/\tilde{p} - 1 + p'/2}$, and

$$
(7.8) \quad L^{(2)}(\varepsilon) \equiv L^{(2)}_{\alpha,p}(\varepsilon) \equiv (1 + \alpha)^{-\gamma_2} \min\left[ 1, \frac{(\alpha + 1)/\log_+(1/\varepsilon)}{1 - 2^{-|(\alpha+1)p'/(p'+1) - 1|}} \right].
$$

THEOREM 7.2.  *Let $R^{(\varepsilon,*)}(\beta)$ be the ideal risk in* (1.3) *and let $R^{(\varepsilon)}(\hat{\beta}^{(\varepsilon)}, \beta)$ be the risk* (1.2) *of the GEB estimator* (2.10). *Then there exists a universal constant $M < \infty$ such that*

$$
\sup \{ R^{(\varepsilon)}(\hat{\beta}^{(\varepsilon)}, \beta) - R^{(\varepsilon,*)}(\beta) : \beta \in B^\alpha_{p,q}(C) \}
$$
(7.9)
$$
\leq MC^2 \left\{ (\varepsilon/C)^2 + \sum_{j=1}^2 (\varepsilon/C)^{2\alpha_j/(\alpha_j + 1/2)} \log_+^{\gamma_j}(C/\varepsilon) L^{(j)}(\varepsilon/C) \right\},
$$

*for all $0 < \varepsilon \leq C$ and Besov balls $B^\alpha_{p,q}(C) \in \mathcal{B}_{\text{Besov}}$ in* (3.8), *where $\alpha_j > \alpha$ and $\gamma_j$ are as in* (7.5) *and* (7.6), *and $L^{(j)}$ are the bounded nonincreasing slowly varying functions in* (7.7) *and* (7.8).

REMARK. (i) Since $\alpha_1 \geq \alpha_2 > \alpha$, the right-hand side of (7.9) is of smaller order than $\varepsilon^{2\alpha/(\alpha+1/2)}$. Thus, (3.10) and (7.9) imply Theorem 3.3.

(ii) The scale equivariance (2.12) and (2.13) of the GEB estimators (2.10) is reflected in (7.9).

7.3. *Minimax risks in Besov balls.* Let $\mathcal{R}^{(\varepsilon)}(B)$ and $R^{(\varepsilon,*)}(\beta)$ be the minimax and ideal risks in (1.4) and (3.1). It follows from Theorem 7.2 that for all Besov balls $B \in \mathcal{B}_{\text{Besov}}$

$$
(7.10) \quad \mathcal{R}^{(\varepsilon)}(B) \leq \sup_{\beta \in B} R^{(\varepsilon,*)}(\beta) + o(1) \varepsilon^{2\alpha/(\alpha + 1/2)} \qquad \text{as } \varepsilon \to 0+.
$$

In this section, we provide an inequality which implies

$$
(7.11) \qquad \lim_{\varepsilon \to 0+} \frac{\sup_{\beta \in B^\alpha_{p,q}(C)} R^{(\varepsilon,*)}(\beta)}{\mathcal{R}^{(\varepsilon)}(B^\alpha_{p,q}(C))} = 1.
$$

Let $\|(x_1,\ldots,x_n)\|_{p,n} \equiv (\sum_{k=1}^n |x_k|^p)^{1/p}$ for $p > 0$ and $n \geq 1$ with usual extensions for $p \vee n = \infty$. Let $C_j$ denote nonnegative constants. It follows from (1.4) and (3.9) that

$$\mathcal{R}^{(\varepsilon)}(B_{p,q}^\alpha(C))$$

$$\geq \sup_{\|\{C_j\}\|_{q,\infty} \leq C} \sum_{j=-1}^\infty \inf_{\hat{\beta}} \sup \left\{ \sum_{k=1}^{2^j \vee 1} E_\beta^{(\varepsilon)}(\hat{\beta}_{jk} - \beta_{jk})^2 : \frac{\|\beta_{[j]}\|_{p,2^j}}{2^{j/p}} \leq \frac{C_j}{2^{j(\alpha+1/2)}} \right\},$$

so that by (6.28) with the scale change $\theta_{jk} = \beta_{jk}/\varepsilon$

$$\mathcal{R}^{(\varepsilon)}(B_{p,q}^\alpha(C)) \geq \varepsilon^2 \sup_{\|\{C_j\}\|_{q,\infty} \leq C} \sum_{j=1}^\infty 2^j \mathcal{R}_{2^j}(\Theta_{p,2^j}(2^{-j(\alpha+1/2)}C_j/\varepsilon)).$$

Furthermore, it follows from (7.1) and (7.4) that

$$\sup_{\beta \in B_{p,q}^\alpha(C)} R^{(\varepsilon,*)}(\beta)$$

$$= \varepsilon^2 \sup_{\|\{C_j\}\|_{q,\infty} \leq C} \sum_{j=1}^\infty 2^j \sup \{R^*(G_{[j]}^{(\varepsilon)}) : \mu_p(G_{[j]}^{(\varepsilon)}) \leq 2^{-j(\alpha+1/2)}C_j/\varepsilon\}.$$

The above facts and Proposition 6.3 imply

$$\sup_{\beta \in B_{p,q}^\alpha(C)} R^{(\varepsilon,*)}(\beta) - \mathcal{R}^{(\varepsilon)}(B_{p,q}^\alpha(C))$$

(7.12)
$$\leq \varepsilon^2 \sup_{\|\{C_j\}\|_{q,\infty} \leq C} \sum_{j=1}^\infty 2^j M r_p^*(2^j, 2^{-j(\alpha+1/2)}C_j/\varepsilon)$$

for the $r_p^*(n,C)$ in (6.31), since $\sup_{\mu_p(G) \leq C} R^*(G) - \mathcal{R}_n(\Theta_{p,n}(C)) \leq M r_p^*(n,C)$ for all $(n,C)$.

Theorem 7.3 below, which implies (7.11), is a consequence of (7.12). Define

(7.13) $\qquad \alpha_3 \equiv \alpha + (\alpha + 1/2)/2, \qquad \gamma_3 \equiv 2/3,$

and for $\gamma \equiv p'(\alpha + 1/2) - 1$ define bounded slowly varying functions $L^{(3)}(\varepsilon) \equiv L_{\alpha,p}^{(3)}(\varepsilon)$ by

(7.14)
$$L^{(3)}(\varepsilon) \equiv \frac{\log_+^{2/3}(\varepsilon^{-1/(\alpha+1/2)}p^2 p')}{(p^2 p')^{1/3} \log_+^{2/3}(1/\varepsilon)}$$
$$\times \left[ (\gamma+1)^{1/3} + \frac{(\gamma+1)^{(2-p')/3}}{(1-2^{-\gamma})^{1+(2-p')/3}} \right.$$
$$\left. + \log_+^{-2/3}(\varepsilon^{-1/(\alpha+1/2)}p^2 p') \frac{(\gamma+1)^{(4-p')/3}}{(1-2^{-\gamma})^{1+(4-p')/3}} \right].$$

THEOREM 7.3. *Let $\mathcal{R}^{(\varepsilon)}(B)$ and $R^{(\varepsilon,*)}(\beta)$ be the minimax and ideal risks in* (1.4) *and* (1.3). *Then* (3.10) *and* (7.11) *hold, and there exists a universal constant $M < \infty$ such that*

$$(7.15) \quad \sup_{\beta \in B} R^{(\varepsilon,*)}(\beta) \leq \mathcal{R}^{(\varepsilon)}(B) + MC^2 (\varepsilon/C)^{2\alpha_3/(\alpha_3+1/2)} \log_+^{\gamma_j}(C/\varepsilon) L^{(3)}(\varepsilon/C)$$

*for all $0 < \varepsilon \leq C$ and Besov balls $B = B_{p,q}^\alpha(C) \in \mathcal{B}_{\text{Besov}}$ in* (3.8), *where $\alpha_3 > \alpha$ and $\gamma_3$ are the constants in* (7.13), *and $L^{(3)}$ is the bounded nonincreasing slowly varying function in* (7.14).

REMARK. For $p \leq q$, Donoho and Johnstone (1998) proved (3.10) and (7.11) using the minimax theorem for certain classes of random $\beta$.

7.4. *Exactly adaptive minimaxity and superefficiency.* The universal exactly adaptive minimaxity and related inequalities in Theorem 7.4 below follow immediately from Theorems 7.2 and 7.3, since the sum of the right-hand sides of (7.9) and (7.15) is of smaller order than the rate $\varepsilon^{2\alpha/(\alpha+1/2)}$ in (3.10), due to $\alpha_j > \alpha$, $j = 1, 2, 3$.

THEOREM 7.4. *Let $\hat{\beta}^{(\varepsilon)}$ be the GEB estimator in* (2.3) *or* (2.10), *and let $\mathcal{R}^{(\varepsilon)}(B)$ be the minimax risk in* (1.4). *Then there exists a universal constant $M < \infty$ such that*

$$\mathcal{R}^{(\varepsilon)}(B_{p,q}^\alpha(C))$$
$$\leq \sup \{ R^{(\varepsilon)}(\hat{\beta}^{(\varepsilon)}, \beta) : \beta \in B_{p,q}^\alpha(C) \}$$
$$\leq \mathcal{R}^{(\varepsilon)}(B_{p,q}^\alpha(C))$$
$$\quad + MC^2 \bigg\{ (\varepsilon/C)^2 + \sum_{j=1}^{3} (\varepsilon/C)^{2\alpha_j/(\alpha_j+1/2)} \log_+^{\gamma_j}(C/\varepsilon) L^{(j)}(\varepsilon/C) \bigg\},$$

*for all $0 < \varepsilon \leq C$ and Besov balls $B_{p,q}^\alpha(C) \in \mathcal{B}_{\text{Besov}}$ in* (3.8), *where constants $\alpha_j > \alpha$ and $\gamma_j$ and bounded functions $L^{(j)}$ are as in Theorems* 7.1 *and* 7.2.

REMARK. Since $\alpha_j > \alpha$, Theorem 7.4 and (3.10) imply Theorem 3.4.

Now we consider the superefficiency of the GEB estimators. Let $B$ be a compact set under the Besov norm $\|\cdot\|_{p,q}^\alpha$ (3.9) with $q < \infty$. Let $\Pi_J$ be the projections up to resolution levels $J$, $(\Pi_J \beta)_{jk} \equiv \beta_{jk} I_{\{j \leq J\}}$. Since $\|\beta - \Pi_J \beta\|_{p,q}^\alpha \to 0$ for every $\beta \in B$ and $B$ is compact,

$$(7.16) \qquad c_J^*(B) \equiv \sup\{\|\beta - \Pi_J \beta\|_{p,q}^\alpha : \beta \in B\} \to 0 \qquad \text{as } J \to \infty.$$

The superefficiency follows, since the risk for the estimation of $\Pi_J\beta$ by the GEB estimator is at most $O(1)2^J\varepsilon^2$ and $\{\beta - \Pi_J\beta : \beta \in B\} \subseteq B_{p,q}^\alpha(c_J^*(B))$. Formally, by (7.1) with $n_j = 2^j$,

$$R^{(\varepsilon,*)}(\beta) \leq \varepsilon^2 2^{J+1} + \varepsilon^2 \sum_{j=J+1}^\infty \sum_{k=1}^{2^j} 2^j R^*(G_{[j]}^{(\varepsilon)}) \leq \varepsilon^2 2^{J+1} + \sup_{\|\beta\|_{p,q}^\alpha \leq c_J^*(B)} R^{(\varepsilon,*)}(\beta)$$

for all $\beta \in B$. Since $c_J(B) \to 0$, the right-hand side above is $o(\varepsilon^{2\alpha/(\alpha+1/2)})$ as $\varepsilon \to 0+$ and then $J \to \infty$, by (7.11) and (3.10) in Theorem 7.3. This and (7.9) imply (3.12) and complete the proof of Theorem 3.5.

**8. Bayes and more general classes.** The results in Sections 3 and 7 can be extended in several directions, for example, Bayes models, more general deterministic and stochastic $\beta$ and blocks with sizes $n_j \neq 2^j$. The extension to stochastic $\beta$ is relatively straightforward, since the key oracle inequalities in Theorems 3.1 and 7.1 are valid under integration over $\beta$, for example,

$$\begin{aligned}
&ER^{(\varepsilon)}(\hat\beta^{(\varepsilon)}, \beta) - ER^{(\varepsilon,*)}(\beta) \\
&\qquad \leq M\varepsilon^2 \sum_j \left\{ n_j r_0(n_j, EG_{[j]}^{(\varepsilon)}) + \frac{1}{(\log n_j + 1)^{3/2}} \right\},
\end{aligned} \tag{8.1}$$

due to the concavity of $r_0(n, G)$ in $G$. We consider here certain general classes of Bayes models including wavelet coefficients in Besov balls and of functions with a large number of discontinuities.

Let $\beta$ be a random sequence and let $E$ be the certain expectation under which $E_\beta^{(\varepsilon)}$ of model (1.1) is the conditional expectation given $\beta$. Let $\vec{\mu}_p(\beta)$ be the sequence $\{(E|\beta_{jk}|^p)^{1/p}\}$ of the marginal $L_p$ norms of $\beta$. Let $\beta_{[j]} \equiv \{\beta_{jk}, k \leq 1 \vee 2^j\}$ and $\kappa(X_1, \ldots, X_n) \equiv n^{-1} \sum_{k=1}^n E(X_k^2 \wedge 1)$. Consider

$$(8.2) \quad B^{(\varepsilon)} \equiv \left\{ \beta = \beta' + \beta'' : \vec{\mu}_p(\beta') \in B_{p,q}^\alpha(C), \kappa(\beta_{[j]}''/\varepsilon) \leq \frac{m^{(\varepsilon)}}{2^j}\left[1 \wedge \frac{M^{(\varepsilon)}}{\varepsilon^2 2^j}\right]\right\},$$

where $m^{(\varepsilon)}$ and $M^{(\varepsilon)}$ are constants. Let $\mathcal{F}_{d,m}(C)$ be the class of piecewise polynomials $f$ of degree $d$ with no more than $m$ pieces and $\|f\|_\infty \leq C$. A deterministic $\beta = \beta' + \beta''$ belongs to $B^{(\varepsilon)}$ if $\beta' \in B_{p,q}^\alpha(C)$ and $\beta_{jk}'' = \int \phi_{jk} f$ are the wavelet coefficients of $f \in \mathcal{F}_{d,cm^{(\varepsilon)}}(cM^{(\varepsilon)})$ as in Example 3.1, for certain fixed small $c > 0$.

THEOREM 8.1. *Suppose $(\log \varepsilon)^2 m^{(\varepsilon)} = o(1)\varepsilon^{-1/(\alpha+1/2)}$ and $\log_+(M^{(\varepsilon)}) = O(|\log \varepsilon|)$. Then (2.10) is uniformly adaptive to the ideal risk $R^{(\varepsilon,*)}(\beta)$ in (1.3) over classes (8.2) of random $\beta$,*

$$(8.3) \quad \sup\{ER^{(\varepsilon)}(\hat\beta^{(\varepsilon)}, \beta) - ER^{(\varepsilon,*)}(\beta) : \beta \in B^{(\varepsilon)}\} = o(\varepsilon^{2\alpha/(\alpha+1/2)}).$$

*Moreover, the GEB estimators are exactly adaptive minimax,*

$$(8.4) \quad \sup_{\beta \in B^{(\varepsilon)}} ER^{(\varepsilon)}(\hat{\beta}^{(\varepsilon)}, \beta) = (1+o(1))\mathcal{R}^{(\varepsilon)}(B^{(\varepsilon)}) = (1+o(1))\mathcal{R}^{(\varepsilon)}(B_{p,q}^{\alpha}(C)),$$

*where $\mathcal{R}^{(\varepsilon)}(B)$ is the minimax $\ell_2$ risk for the estimation of a random $\beta$ in $B$.*

REMARK. (i) Although $B^{(\varepsilon)}$ in (8.2) is much larger than the Besov class $B_{p,q}^{\alpha}(C)$, the minimax risks for the two classes are within an infinitesimal fraction of each other.

(ii) The condition on $m^{(\varepsilon)}$ in Theorem 8.1 is the weakest possible up to a factor of $(\log \varepsilon)^2$, since $m^{(\varepsilon)} = o(1)\varepsilon^{-1/(\alpha+1/2)}$ is a necessary condition.

(iii) Deterministic versions of the classes (8.2) were considered in Hall, Kerkyacharian and Picard (1998) in the context of density estimation.

**9. Equivalence between the white noise model and nonparametric regression.** In this section we establish the asymptotic equivalence between the problems of estimating $f$ in the nonparametric regression model (4.1) and $\beta$ in the white noise model (1.1) in Besov classes when $\beta \equiv \beta(f)$ are the Haar coefficients of $f$. The asymptotic equivalence is used to prove the adaptive minimaxity of the GEB estimators in Theorem 4.1. We assume throughout the section that the design variables $t_i$ in (4.1) are independent uniformly distributed in $(0, 1)$.

THEOREM 9.1. *Let $E_f$ be as in Theorem 4.1 with i.i.d. uniform $t_i$. Let $\beta(f)$ and $\tilde{\beta}(f)$ be as in (4.5) and (4.7), respectively, and let $\Pi_J: \beta_{j,k} \to \beta_{j,k}I\{j \leq J\}$ be the projections as in (7.16).*

(i) *There exist finite constants $M_{\alpha,p}$ such that*

$$(9.1) \quad E_f\{\|\Pi_J\tilde{\beta}(f) - \Pi_J\beta(f)\|_{p',p'}^{\alpha}\}^{p'} \leq M_{\alpha,p'}(2^J/N)^{p'/2}\{\|\beta(f)\|_{p,q}^{\alpha}\}^{p'},$$

*where $p' = p \wedge 2$, and*

$$(9.2) \quad \begin{aligned} &E_f\|\Pi_J\tilde{\beta} - \beta\|_2^2 \\ &\leq M_{\alpha,p'}(\|\beta\|_{p,q}^{\alpha})^2\left\{\frac{J}{N}I\{\alpha p'=1\} + \frac{1}{N} + \frac{2^J/N+1}{2^{2J(\alpha+1/2-1/p')}}\right\}. \end{aligned}$$

(ii) *Let $\varepsilon = \sigma/\sqrt{N}$ and $N \to \infty$. For $\mathcal{F} \equiv \{f: \beta(f) \in B_{p,q}^{\alpha}(C)\}$ and estimates $\hat{f}_N$ based on (4.1),*

$$(9.3) \quad \inf_{\hat{f}_N} \sup_{f \in \mathcal{F}} E_f\|\hat{f}_N - f\|^2 = (1+o(1))\mathcal{R}^{(\varepsilon)}(B_{p,q}^{\alpha}(C))$$

*for $\alpha^2/(\alpha+1/2) > 1/p' - 1/2$, where $\|f\| \equiv (\int_0^1 f^2)^{1/2}$.*

Theorem 9.1(i) provides upper bounds for the difference between the wavelet coefficients $\beta_{j,k} = \int f\phi_{jk}$ and the corresponding coefficients $\tilde{\beta}_{j,k}$ for the random discrete Haar system in (4.7). Deterministic discrete wavelet systems were considered in Donoho and Johnstone (1998) based on Dubuc (1986). For $\alpha > 1/p$ and $p \vee q \geq 1$ or $\alpha = p = q = 1$, Donoho and Johnstone (1998) established (9.3) for deterministic discrete wavelet systems.

## APPENDIX

We shall denote by $M$ generic finite universal constants which may take different values from one appearance to the next, that is, $M \equiv O(1)$ uniformly.

PROOF OF THEOREM 3.2. Consider small $\varepsilon > 0$. Let $\eta > 0$ and

$$\tilde{r}_p(n, C) \equiv \min[1, C^p, (C/\sqrt{n})^{p/(p+1)}]$$
$$= \min[1, C^p, \max\{n^{-1/2}, (C/\sqrt{n})^{p/(p+1)}\}].$$

It follows from Theorem 3.1 and part three of (3.4) that the regret (3.6) is bounded by

$$\sup_{\|\beta\| \leq C} r^{(\varepsilon)}(\hat{\beta}^{(\varepsilon)}, \beta) \leq O(\varepsilon^2) + o(\varepsilon^{2-\eta}) \sum_j n_j \tilde{r}_{p'}(n_j, (\varepsilon n_j^s)^{-1}), \qquad p' \equiv p \wedge 2.$$

We compute the above bound by splitting the sum into three pieces for $n_j \in [x_k, x_{k+1})$, $k = 0, 1, 2$, where $x_0 = 1$, $x_1 = \varepsilon^{-1/(s+1/2)}$, $x_2 = \varepsilon^{-2/(2s-1/p')}$ and $x_3 = \infty$. This yields by (3.4)

$$\sum_j n_j \tilde{r}_{p'}(n_j, (\varepsilon n_j^s)^{-1})$$

$$\leq \sum_{j < x_1} n_j + \sum_{x_1 \leq n_j < x_2} n_j (\varepsilon n_j^{s+1/2})^{-p'/(p'+1)} + \sum_{x_2 \leq n_j} n_j (\varepsilon n_j^s)^{-p'}$$

$$\leq o(\varepsilon^{-\eta})\{x_1 + x_1 + x_2(\varepsilon x_2^s)^{-p'}\} = o(\varepsilon^{-\eta - v}),$$

where $v = \max\{1/(s+1/2), 1/(2s - 1/p')\} = 1/(\alpha_0 + 1/2)$. Thus the regret is uniformly bounded by $o(1)\varepsilon^{2-v-2\eta}$. This completes the proof, since $\eta$ is arbitrary and $2 - v = 2\alpha_0/(\alpha_0 + 1/2)$. $\square$

LEMMA A.1. *Let $h^{(m)}(x) = (d/dx)^m h(x)$. Let $\hat{\varphi}_n$ be given by (6.5) with $a \equiv a_n \geq \sqrt{2\log n}$. For $p \geq 2$, there exist universal constants $M_p < \infty$ such that*

$$\int E|\hat{\varphi}_n^{(m)}(x) - \varphi_G^{(m)}(x)|^p dx \leq M_p \frac{a^{mp+p/2}}{n^{p/2}} \left\{1 + \left(\frac{a}{\sqrt{n}}\right)^{p/2-1}\right\}.$$

PROOF. Let $M_p$ denote any positive universal constant. We shall omit the calculation involving the bias $b_n(x) \equiv E\hat{\varphi}_n(x) - \varphi_{G_n}(x)$, since it is of smaller order in the sense that

$$\|b_n^{(m)}(x)\|_\infty \leq \frac{1}{\pi} \int_a^\infty u^m e^{-u^2/2}\, du \leq O(1)\frac{a^{m-1}}{n}$$

by (6.5) and the Fourier inversion formula, and by the Plancherel identity

$$\int |b_n^{(m)}(x)|^2\, dx = \frac{1}{\pi}\int_a^\infty u^{2m} e^{-u^2}\, du \leq O(1)\frac{a^{2m-1}}{n^2}.$$

Let $W_k(x) \equiv a^{m+1}K^{(m)}(a(x-X_k))$ and $h_p(x) \equiv \sum_{k=1}^n E|W_k(x)|^p/n$. Since $\hat{\varphi}_n^{(m)}(x)$ is the average of $W_k(x)$ and $\{W_k(x), k \leq n\}$ are independent given $\{\theta_k, k \leq n\}$,

$$E|\hat{\varphi}_n^{(m)}(x) - E\hat{\varphi}_n^{(m)}(x)|^p \leq \frac{M_p}{n^p} h_2^{p/2}(x) + \frac{M_p}{n^p} h_p(x).$$

This implies the conclusion, since $\|h_p(x)\|_\infty + \int h_p(x)\, dx = O(a^{mp+p-1})$ via

$$h_p(x) = \int |a^{m+1}K^{(m)}(a(x-u))|^p \varphi_{G_n}(u)\, du$$

$$= a^{mp+p-1}\int |K^{(m)}(u)|^p \varphi_{G_n}(x-u/a)\, du. \qquad \square$$

PROOF OF THEOREM 6.1. The difference between the proof here and that of Zhang (1997) is the use of the improved bounds in Lemma A.1. We shall only describe the differences and refer to Zhang (1997) for the rest. Let $a \equiv a_n = \sqrt{2\log n}$.

The condition $\rho^{-1}a/\sqrt{n} = o(1)$ of Lemma 1 of Zhang (1997) can be weakened to $\rho^{-1}\sqrt{a/n} = o(1)$, since by Theorem 2 of Zhang (1997) and Lemma A.1

$$E\int \{\hat{\varphi}_n^{(m)}(x) - \varphi_G^{(m)}(x)\}^2 \frac{\max(\varphi_{G_n}(x), \rho)}{\max(\hat{\varphi}_n(x), \rho)}\, dx$$

$$\leq E\int \{\hat{\varphi}_n^{(m)}(x) - \varphi_G^{(m)}(x)\}^2\{1 + |\hat{\varphi}_n(x) - \varphi_{G_n}(x)|/\rho\}\, dx$$

$$\leq E\|\hat{\varphi}_n^{(m)} - \varphi_G^{(m)}\|_2^2 + \rho^{-1}\sqrt{E\|\hat{\varphi}_n^{(m)} - \varphi_G^{(m)}\|_4^4 E\|\hat{\varphi}_n - \varphi_{G_n}\|_2^2}$$

$$\leq \frac{(1+o(1))a^{2m+1}}{(2m+1)\pi n} + \frac{Ma^{2m+3/2}}{\rho n^{3/2}}.$$

The assumption $\rho^{-1}a/\sqrt{n}$ is used in the proof of Theorem 3 of Zhang (1997) only for the application of Lemma 1 there. The assumption $a = O(\sqrt{\log n})$ is actually not used in the proof of Theorem 3 of Zhang (1997).

Hence, Theorem 3 of Zhang (1997) holds under weaker conditions $a \geq \sqrt{\log n}$ and $\rho^{-1}\sqrt{a/n} = o(1)$. The proof in Section 5 of Zhang (1997) is based on Theorem 3 of Zhang (1997) and the additional conditions $a \geq \sqrt{2\log n}$ and $\Delta^*(n,\rho) = O(1)$ only, since $a^{3/2}/(\rho n) = o(1) a^3/(\rho n)$. □

PROOF OF THEOREM 6.2. Part (i) follows directly from Remark 2 below Theorem 6.1 and the fact that $\Delta^*(n,\rho_n) \to 0$ in Theorem 6.1. For part (ii), we shall first prove (6.11). Let $\varphi_j(\cdot) \equiv \varphi(\cdot; H_j)$ and $\varphi_G \equiv \varphi(\cdot; G)$ be as in (2.5) and $\tilde{w}_j \equiv w_j \varphi_j / \varphi_G$. Since $\sum_j \tilde{w}_j = 1$, by Cauchy–Schwarz

$$\left(\frac{\varphi'_G}{\varphi_G}\right)^2 \varphi_G = \left(\sum_{j=0}^m \tilde{w}_j \frac{\varphi'_j}{\varphi_j}\right)^2 \varphi_G \leq \sum_{j=0}^m \tilde{w}_j \left(\frac{\varphi'_j}{\varphi_j}\right)^2 \varphi_G = \sum_{j=0}^m w_j \left(\frac{\varphi'_j}{\varphi_j}\right)^2 \varphi_j.$$

This and (6.7) imply (6.11), since $1 - \varphi_G/(\varphi_G \vee \rho)$ is decreasing in $\varphi_G$ and $\varphi_G \geq w_j \varphi_j$. Let $A$ be a union of $m$ disjoint intervals $I_j$ of length $\leq 1$. A distribution $G$ can be written as $G = \sum_{j=0}^m w_j H_j$, where $w_0 = G(A^c)$, $w_j \equiv G(I_j)$ and $H_j$ are the conditional distributions given $\theta \in I_j$ under $G$. Define $\eta(\rho) \equiv \sup\{\Delta(\rho, H) : H([0,1]) = 1\}$. Since $\Delta(\rho, H_j) \leq 1$, (6.11) implies

$$\Delta(\rho, G) \leq G(A^c) + \sum_{j=1}^m w_j I_{\{\rho/w_j \geq 1/M\}} + \eta(1/M)$$

(A.1)
$$\leq G(A^c) + Mm\rho + \eta(1/M).$$

It follows from Proposition 2 of Zhang (1997) that $\eta(\rho) \to 0$ as $\rho \to 0$.

Now, let $A_n \equiv \bigcup_{j=1}^{m_n}[c_{n,j} - M, c_{n,j} + M]$. The condition of part (ii) implies $G_n(A_n^c) \leq w_{n,0} + \sup\{H([-M,M]^c) : H \in \mathcal{G}\} \to 0$ for large $n$ and $M$. Thus, we may assume $G_n(A_n^c) \to 0$ for certain $A_n = \bigcup_{j=1}^{m_n} I_{n,j}$ with disjoint intervals $\{I_{n,j}, j \leq m_n\}$ of at most unit length and (possibly different) $m_n = o(1/\rho_n)$. Under this assumption and conditionally on $\theta_{(n)}$, Theorem 6.1 and (A.1) imply

$$\frac{1}{n}\sum_{k=1}^n E(\hat{\theta}_k - \theta_k)^2 \leq ER^*(G_{(n)}) + EG_{(n)}(A_n^c) + Mm_n\rho_n + \eta(1/M) + o(1)$$

$$\leq ER^*(G_n) + o(1)$$

with $G_{(n)}(x) \equiv n^{-1}\sum_{k=1}^n I\{\theta_k \leq x\}$, as $n \to \infty$ and then $M \to \infty$, since $ER^*(G_{(n)}) \leq R^*(EG_{(n)}) = R^*(G_n)$ due to the concavity of $R^*(G)$ in $G$ and $EG_{(n)}(A_n^c) = G_n(A_n^c) \to 0$. □

PROOF OF LEMMA 6.1. Let $x$ be fixed. Let $H_1$ and $H_2$ be the conditional distributions given $|\theta| > x$ and $|\theta| \leq x$, respectively, under $G$. Let $w_1 \equiv \overline{G}(x)$, $w_2 \equiv 1 - w_1$ and $\varphi_j(\cdot) \equiv \varphi(\cdot; H_j)$ be as in (2.5). Since $H_2([-x,x]) = 1$, by the unimodality of $\varphi$, $\varphi_2$ is monotone in both $(-\infty, -x)$ and $(x, \infty)$. By

Lemma 2 of Zhang (1997), $|\varphi_2'/\varphi_2| \leq \widetilde{L}(\varphi_2)$. This and the monotonicity of $\varphi_2$ imply

$$\int_x^\infty I_{\{\varphi_2<\rho/w_2\}}\left(\frac{\varphi_2'}{\varphi_2}\right)^2 \varphi_2 \leq \int_x^\infty I_{\{\varphi_2<\rho/w_2\}}\widetilde{L}(\varphi_2)|\,d\varphi_2|$$

$$\leq \int_0^{\rho/w_2} \widetilde{L}(u)\,du \leq \frac{1}{w_2}\int_0^\rho \widetilde{L}(u)\,du$$

and a similar inequality for $\int_{-\infty}^{-x}$. These and Cauchy–Schwarz imply

$$\frac{w_2}{2}\int_{|u|>x} I_{\{\varphi_2<\rho/w_2\}}\left(\frac{\varphi_2'}{\varphi_2}\right)^2 \varphi_2\,du \leq \int_0^\rho \widetilde{L}(u)\,du \leq \left(\rho\int_0^\rho \widetilde{L}^2(u)\,du\right)^{1/2}$$

(A.2)
$$= \rho\sqrt{\widetilde{L}^2(\rho)+2}.$$

For (6.13), we find again by Lemma 2 of Zhang (1997) that

$$w_2\int_{|u|\leq x} I_{\{\varphi_2<\rho/w_2\}}\left(\frac{\varphi_2'}{\varphi_2}\right)^2 \varphi_2\,du \leq w_2\int_{|u|\leq x} I_{\{\varphi_2<\rho/w_2\}}\widetilde{L}^2(\varphi_2)\varphi_2\,du$$

$$\leq 2x\max\{\widetilde{L}^2(\rho),2\}\rho$$

due to the monotonicity of $u\max\{\widetilde{L}^2(u),2\}$ and $\widetilde{L}^2(u)$. Thus, (6.13) holds, since by (6.11)

$$\Delta(\rho,G) \leq w_1 + w_2\Delta(\rho/w_2,H_2) \leq w_1 + w_2\int_{\varphi_2\leq\rho/w_2}\left(\frac{\varphi_2'}{\varphi_2}\right)^2 \varphi_2.$$

For (6.14), $\{1-\varphi_2/(\varphi_2\vee(\rho/w_2))\} \leq \{1-\rho/(\rho\vee(\rho/w_2))\} = w_1$ for $|u|\leq x$, since $\varphi_2(u)\geq\varphi(2x)=\rho$. Thus, (6.14) follows from (A.2) and

$$w_2\int_{|u|\leq x}\left(\frac{\varphi_2'}{\varphi_2}\right)^2 \varphi_2\left(1-\frac{\varphi_2}{\varphi_2\vee(\rho/w_2)}\right)^2 \leq w_1^2 w_2. \qquad \square$$

PROOF OF THEOREM 6.3. Since the right-hand sides of (6.16) and (6.17) are both concave in $G$, it suffices to apply Theorem 6.1 conditionally on $\theta_{(n)}$. By (6.8) and simple calculation, (6.15) holds, so that (6.17) follows from Theorem 6.1 and (6.14). For (6.16) we use the Markov inequality $\overline{G}(x) \leq \mu_p^p(G)/x^p$ in (6.13) and then minimize $\mu_p^p(G)/x^p + 2x\rho_n(1+o(1))\log n$ over $x>0$. $\square$

PROOF OF LEMMA 6.2. By (6.1) it suffices to verify (6.19) for degenerate $G$. Let $X \sim N(\mu,1)$ under $P_\mu$. Let $h(x)\equiv s_\lambda(x)-x = \min(\lambda,\max(-\lambda,-x))$.

For $\theta_1 = \mu$,
$$R(s_\lambda, G_1) = R_s(\mu; \lambda) \equiv E_\mu\{h(X) + X - \mu\}^2 = E_0\{h(X + \mu) + X\}^2.$$

Differentiating twice the right-hand side above with respect to $\mu$, we find
$$\left(\frac{\partial}{\partial \mu}\right)^2 R_s(\mu; \lambda) = 2\left[\int_{-\lambda-\mu}^{\lambda-\mu} \varphi(u)\,du + \mu\varphi(\lambda+\mu) - \mu\varphi(\lambda-\mu)\right] \leq 2$$

for all positive $\mu$ and $\lambda$. Since $R_s(\mu; \lambda)$ is an even function, $R_s(\mu; \lambda) \leq R_s(0; \lambda) + \mu^2$. This implies the first component of (6.19) due to
$$R_s(0; \lambda) = 2\int_\lambda^\infty (u-\lambda)^2 \varphi(u)\,du$$
$$= 2\varphi(\lambda) \int_0^\infty u^2 e^{-\lambda u - u^2/2}\,du \leq 2\lambda^{-3}\varphi(\lambda) \int_0^\infty u^2 e^{-u}\,du.$$

The second component of (6.19) follows from the monotonicity of $R_s(\mu; \lambda)$ in $|\mu|$, proved below, as $\lim_{\mu \to \infty} R_s(\mu; \lambda) = \lambda^2 + 1$. By Stein's formula of mean-squared error,
$$R_s(\mu; \lambda) = E_\mu\{h^2(X) + 1 + 2h'(X)\}$$
$$= \int_0^\lambda P_\mu\{|X| > u\}\,du^2 + 2P_\mu\{|X| > \lambda\} - 1.$$

The monotonicity of $R_s(\mu; \lambda)$ then follows from that of $P_\mu\{|X| > u\}$ in $|\mu|$. Inequality (6.20) is a direct consequence of (6.19). □

LEMMA A.2. *Let $U_k$ be independent random variables with $P\{0 \leq U_k \leq 1\} = 1$. Set $\mu_n \equiv n^{-1} \sum_{k=1}^n EU_k$. For all $0 < \mu_n < u < 1$,*
$$\text{(A.3)} \quad P\left\{n^{-1} \sum_{k=1}^n U_k > u\right\} \leq \exp[-nK(u, \mu_n)] \leq \exp[-2n(u-\mu_n)^2],$$

*where $K(p_1, p_2)$ is the Kullback–Leibler information for Bernoulli variables, defined by*
$$K(p_1, p_2) \equiv p_1 \log\left(\frac{p_1}{p_2}\right) + (1-p_1)\log\left(\frac{1-p_1}{1-p_2}\right)$$
$$= \int_0^{p_1-p_2} \frac{p_1 - p_2 - u}{(p_2+u)(1-p_2-u)}\,du.$$

PROOF. Let $p_k \equiv EU_k$ and $\delta_k$ be Bernoulli variables with $E\delta_k = p_k$. Since $EU_k^m \leq p_k = E\delta_k^m$ for all integer $m \geq 0$ and $\log(1 + p_k(e^\lambda - 1))$ is concave in $p_k$, for $\lambda \geq 0$,

$$E\exp\left(\lambda\sum_{k=1}^{n}U_k\right)\leq\prod_{k=1}^{n}Ee^{\lambda\delta_k}=\prod_{k=1}^{n}(1+p_k(e^\lambda-1))\leq(1+\mu_n(e^\lambda-1))^n.$$

The first inequality of (A.3) follows from $P\{\sum_k U_k > nu\} \leq e^{-\lambda nu}(1+\mu_n(e^\lambda-1))^n$ with $\lambda = \log[\{u(1-\mu_n)\}/\{\mu_n(1-u)\}]$. The second one follows from the integral formula of the Kullback–Leibler information and the bound $(p_2+u)(1-p_2-u) \leq 1/4$. □

PROOF OF LEMMA 6.3. By (6.21) and (6.23), $E\hat{\kappa}_n = \tilde{\kappa}(G_n)$, so that by Lemma A.2,

$$(A.4) \qquad P\{\pm(\hat{\kappa}_n - \tilde{\kappa}(G_n)) > u\} \leq \exp(-nu^2) \qquad \forall u > 0,$$

with $U_k$ (or $1-U_k$) being $\exp(-X_k^2/2)$. Since $\delta_n \equiv I\{\hat{\kappa}_n \leq b_n\}$ are Bernoulli variables,

$$\frac{1}{n}\sum_{k=1}^{n}E(\hat{t}_n(X_k)-\theta_k)^2$$

(A.5)

$$= \frac{1}{n}\sum_{k=1}^{n}E(\hat{t}_{n,\rho}(X_k)-\theta_k)^2(1-\delta_n) + \frac{1}{n}\sum_{k=1}^{n}E(s_\lambda(X_k)-\theta_k)^2\delta_n.$$

Thus, it suffices to consider the first two cases of (6.25).

Suppose $\tilde{\kappa}(G_n) > b_n^+ \equiv b_n + \sqrt{2(\log n)/n}$. By (A.4)

$$P\{\delta_n = 1\} \leq P\{\hat{\kappa}_n - \tilde{\kappa}(G_n) \leq -\sqrt{2(\log n)/n}\} \leq \exp(-2\log n) = n^{-2},$$

so that by (6.18), with $\chi_n^2 \equiv \sum_{k=1}^n (X_k - \theta_k)^2$ and $\lambda = \sqrt{2\log n}$,

$$E\delta_n \sum_{k=1}^{n}(s_\lambda(X_k) - \theta_k)^2/n \leq E\delta_n(\sqrt{\chi_n^2/n} + \lambda)^2$$

$$\leq \int_{u_{2,n}}^{\infty}(\sqrt{u/n}+\lambda)^2 p_n(u)\,du,$$

where $p_n(u) \equiv (u/2)^{n/2-1}e^{-u/2}/\{2\Gamma(n/2)\}$ is the density of $\chi_n^2$ and $P\{\chi_n^2 > u_{j,n}\} = 1/n^j$. By standard large deviation theory, $u_{j,n} = n + (2+o(1)) \times \sqrt{jn\log n}$ for each $j$. Integration by parts yields $\int_{u_{2,n}}^{\infty}(u/n)p_n(u)\,du \leq (u_{j,n}/n + 1)\int_{u_{2,n}}^{\infty}p_n(u)\,du = (2+o(1))/n^j$. Thus,

$$E\delta_n \sum_{k=1}^{n}(s_\lambda(X_k)-\theta_k)^2/n \leq n^{-2}(\lambda+O(1))^2 = (2+o(1))(\log n)/n^2.$$

Now consider the case $\tilde{\kappa}(G_n) \leq b_n^- \equiv b_n - \sqrt{3(\log n)/n}$. By Lemma A.2

$$P\{\delta_n = 0\} \leq P\{\hat{\kappa}_n - \tilde{\kappa}(G_n) \geq \sqrt{3(\log n)/n}\} \leq n^{-3}.$$

By (6.5) and (2.11) and $a_n = \sqrt{2\log n}$, $|\hat{t}_{n,\rho}(X_k) - X_k| \leq a_n^2/(2\pi\rho) = (\log n)/(\pi\rho)$, so that

$$E(1-\delta_n)\sum_{k=1}^n (\hat{t}_{n,\rho}(X_k) - \theta_k)^2/n$$
$$\leq \int_{u_{3,n}}^\infty (\sqrt{u/n} + (\log n)/(\pi\rho))^2 p_n(u)\, du$$
$$= n^{-3}((\log n)/(\pi\rho) + O(1))^2. \qquad \Box$$

PROOF OF PROPOSITION 6.1. Part (i) follows from the proof of Theorem 6.3. For part (ii) we have

$$\int_0^{\log n} w'\overline{G}'(\sqrt{u})\, du + \inf_{x\geq 1}\left[w''\overline{G}''(x) + x\frac{(\log n)^{3/2}}{\sqrt{n}}\right]$$
$$\geq \inf_{x\geq 1}\left[\overline{G}(x) + x\frac{(\log n)^{3/2}}{\sqrt{n}}\right], \qquad n > 2. \qquad \Box$$

PROOF OF THEOREM 6.4. By Proposition 6.1(i) it suffices to consider $r_0(n,G)$ in the minimum in (6.26) and independent $\{\theta_k\}$. By Lemma 6.3 it suffices to bound $R(s_\lambda, G_n)$ for $\tilde{\kappa}(G_n) < b_n^+$ and $\sum_{k=1}^n E(\hat{t}_{n,\rho}(X_k) - \theta_k)^2/n - R^*(G_n)$ for $\tilde{\kappa}(G_n) > b_n^-$. In fact, by Lemma 6.2 we need

$$(\text{A.6}) \qquad R(s_\lambda, G_n) \leq M\frac{(\log n)^2}{\sqrt{n}}, \qquad \tilde{\kappa}(G_n) < b_n^+,$$

and by Theorems 6.2(i) and 6.3 and the fact that $\kappa(G) = \int_0^1 \overline{G}(\sqrt{u})\, du$ we need

$$(\text{A.7}) \quad \frac{1}{n}\sum_{k=1}^n E(\hat{t}_{n,\rho}(X_k) - \theta_k)^2 - R^*(G_n) \leq M\kappa(G_n), \qquad \tilde{\kappa}(G_n) > b_n^-.$$

By (6.22) and the second part of (6.20), $\tilde{\kappa}(G_n) < b_n^+ = (2 + o(1))(\log n)/\sqrt{n}$ implies

$$R(s_\lambda, G_n) \leq (2\log n + 1)\kappa(G_n) + \frac{1}{n} \leq M\frac{(\log n)^2}{\sqrt{n}},$$

so that (A.6) holds. By (6.12) and (2.7) $\overline{G}_n(1) \leq \kappa(G_n)$, so that

$$\frac{1}{n}\sum_{k=1}^n E(\hat{t}_{n,\rho}(X_k) - \theta_k)^2 - R^*(G_n) \leq \frac{5}{4}\kappa(G_n) + O(1)b_n^-$$

by (6.17) and the fact that $b_n^- \sim (\log n)/\sqrt{n}$. This implies (A.7), since $\tilde{\kappa}(G_n) \leq \kappa(G_n)$ by (6.22). $\Box$

PROOF OF PROPOSITION 6.2. Let $\varphi_k \equiv \varphi(\cdot; H_k)$, $k = 1, 2$, $\varphi_G \equiv (1 - w)\varphi_1 + w\varphi_2$, and $G \equiv (1 - w)H_1 + wH_2$. By (6.9) and algebra

$$R^*(G) - R^*(H_1)$$

$$= \int \left(\frac{\varphi_1'}{\varphi_1}\right)^2 \varphi_1 - \int \left(\sum_{k=1}^2 \tilde{w}_k \frac{\varphi_k'}{\varphi_k}\right)^2 \varphi_G$$

$$= w \int \left(\frac{\varphi_1'}{\varphi_1}\right)^2 \varphi_1 + \int \left(\frac{\varphi_1'}{\varphi_1}\right)^2 \tilde{w}_1 \tilde{w}_2 \varphi_G$$

$$- 2 \int \left(\frac{\varphi_1'}{\varphi_1}\right) \left(\frac{\varphi_2'}{\varphi_2}\right) \tilde{w}_1 \tilde{w}_2 \varphi_G - \int \left(\frac{\varphi_2'}{\varphi_2}\right)^2 \tilde{w}_2^2 \varphi_G,$$

where $\tilde{w}_1 \equiv (1-w)\varphi_1/\varphi_G \in [0,1]$ and $\tilde{w}_2 \equiv 1 - \tilde{w}_1 = w\varphi_2/\varphi_G$. Set $q = \log(\sqrt{2}/w)$. For $q \leq 1$ the right-hand side of (6.29) is greater than $\sup_G R^*(G) = 1$. Assume $q > 1$. By Hölder

$$\int \left(\frac{\varphi_k'}{\varphi_k}\right)^2 \tilde{w}_k \tilde{w}_2 \varphi_G \leq \left[\int \left(\frac{\varphi_k'}{\varphi_k}\right)^{2q} \tilde{w}_k \varphi_G\right]^{1/q} \left[\int \tilde{w}_2 \varphi_G\right]^{1-1/q}$$

$$\leq (E|Z|^{2q})^{1/q} w^{1-1/q}$$

with $Z \sim N(0,1)$, since $\varphi_k'/\varphi_k$ is the conditional expectation of $Z$ given a random variable with density $\varphi_k$. Similarly, due to $\int (\varphi_k'/\varphi_k)^2 \varphi_k \leq 1$,

$$\int \left(\frac{\varphi_1'}{\varphi_1}\right) \left(\frac{\varphi_2'}{\varphi_2}\right) \tilde{w}_1 \tilde{w}_2 \varphi_G \leq \left[\int \left(\frac{\varphi_1'}{\varphi_1}\right)^2 \tilde{w}_1 \tilde{w}_2 \varphi_G \int \left(\frac{\varphi_2'}{\varphi_2}\right)^2 w \varphi_2\right]^{1/2}$$

$$\leq (E|Z|^{2q})^{1/(2q)} w^{1-1/(2q)}.$$

Thus, $|R^*(G) - R^*(H_1)| \leq w(1 + w^{-1/(2q)} \|Z\|_{2q})^2$. Let $h_0(q) \equiv \Gamma(q+1/2)e^q/q^q$. Since $h_0(q) \leq h_0(q+1) \to \sqrt{2\pi}$, $\|Z\|_{2q}^{2q} = \Gamma(q+1/2)2^q/\sqrt{\pi} \leq \sqrt{2}(2q/e)^q$. These two inequalities imply

$$|R^*(G) - R^*(H_1)| \leq w\{1 + (\sqrt{2}/w)^{1/(2q)} \sqrt{2q/e}\}^2$$

$$= w\{1 + \sqrt{2\log(\sqrt{2}/w)}\}^2.$$

Now we prove (6.30). Let $U \equiv \tilde{\theta}_1 - \tilde{\theta}_0$ and $H_t$ be the conditional distribution of $\tilde{\theta}_t \equiv (1-t)\tilde{\theta}_0 + t\tilde{\theta}_1$ given $|U| \leq \eta_2$. For $k = 0, 1$, $G_k$ are mixtures of $H_k$ and the conditional distributions of $\tilde{\theta}_k$ given $|U| > \eta_2$, so that $|R^*(G_k) - R^*(H_k)|$ are bounded by the right-hand side of (6.29) with $w = \eta_1$. Thus, (6.30) follows from $|(d/dt)R^*(H_t)| \leq \sqrt{8}\{1 + 1/\sqrt{\pi}\}\eta_2$. By (6.9) and calculus,

$$(d/dt)R^*(H_t) = -(d/dt)\int [\{E_*(x-\tilde{\theta}_t)\varphi(x-\tilde{\theta}_t)\}^2/E_*\varphi(x-\tilde{\theta}_t)]\,dx$$

$$= E_*[2E_{*,t}Z\{E_{*,t}U(1-Z^2)\} - \{E_{*,t}Z\}^2\{E_{*,t}UZ\}],$$

where $E_*$ is the conditional expectation given $|U| \leq \eta_2$, $Z$ is an $N(0,1)$ variable independent of $(\tilde{\theta}_0, \tilde{\theta}_1)$ and $E_{*,t}$ is the conditional expectation given $Z + \tilde{\theta}_t$ and $|U| \leq \eta_2$. Hence,

$$|(d/dt)R^*(H_t)| \leq \eta_2\{2\sqrt{E_*(1-Z^2)^2} + E_*|Z|^3\} = \eta_2\sqrt{8}\{1 + 1/\sqrt{\pi}\,\}. \quad \square$$

In addition to Proposition 6.2, we need the following lemma for the proof of Proposition 6.3.

LEMMA A.3. *For $p = \infty$, $\sup_{\mu_p(G) \leq C} R^*(G) = \mathcal{R}_n(\Theta_{p,n}(C))$. For $0 < p < \infty$,*

$$\sup_{\mu_p(G) \leq C/b} R^*(G) - \mathcal{R}_n(\Theta_{p,n}(C))$$

$$\leq 2\pi_0\{1 + \sqrt{2\log(\sqrt{2}/\pi_0)}\,\}^2 + \frac{I\{b^p\pi_0 < 1\}4C^2/(b\pi_0^{1/p})^2}{\exp[nK(b^p\pi_0, \pi_0)]},$$

*for all $b > 1$ and $0 < \pi_0 < 1$, where $K(p_1, p_2)$ is the Kullback–Leibler information in Lemma A.2.*

PROOF. Let $\theta \sim G$ with $\mu_p(G) \leq C/b$. Let $G_1$ be the distribution of $\tilde{\theta} \equiv \theta I\{|\theta| \leq M\}$, where $M \equiv C/(b\pi_0^{1/p})$. Since $P\{|\theta - \tilde{\theta}| > 0\} = P\{|\theta| > M\} \leq \mu_p^p(G)/M^p \leq \pi_0$, by Proposition 6.2

$$(A.8) \qquad R^*(G) \leq R^*(G_1) + 2\pi_0\{1 + \sqrt{2\log(\sqrt{2}/\pi_0)}\,\}^2.$$

Let $\nu_n$ be the prior in $\mathbb{R}^n$ under which $\theta_k$ are i.i.d. variables with marginal distribution $G_1$. For $b > 1$ and estimators $\hat{\theta}_{(n)} \in \Theta_{\infty,n}(M)$,

$$E_{\nu_n} \frac{1}{n} E_{\theta_{(n)}} \|\hat{\theta}_{(n)} - \theta_{(n)}\|^2_{2,n}$$

$$\leq \sup\left\{\frac{1}{n} E_{\theta_{(n)}} \|\hat{\theta}_{(n)} - \theta_{(n)}\|^2_{2,n} : \theta_{(n)} \in \Theta_{\infty,n}(M) \cap \Theta_{p,n}(C)\right\}$$

$$+ 4M^2 \nu_n\{n^{-1/p}\|\theta_{(n)}\|_{n,p} > C\}.$$

Taking the infimum on both sides above over $\hat{\theta}_{(n)} \in \Theta_{\infty,n}(M)$, we find by (6.1) that

$$(A.9) \quad R^*(G_1) \leq \mathcal{R}_n(\Theta_{\infty,n}(M) \cap \Theta_{p,n}(C)) + \frac{4C^2}{b^2\pi_0^{2/p}}\nu_n\left\{\sum_{k=1}^n \frac{|\theta_k|^p}{n} > C^p\right\},$$

since all admissible estimators are almost surely in $\Theta_{\infty,n}(M)$ when $\Theta_{\infty,n}(M)$ is the parameter space. Since $|\theta_k|^p/M^p \leq 1$ are i.i.d. variables under $\nu_n$ with

$E_{\nu_n}|\theta_k|^p/M^p \le \pi_0$, by Lemma A.2

$$\text{(A.10)} \quad \nu_n\left\{\sum_{k=1}^n \frac{|\theta_k|^p}{M^p n} > \frac{C^p}{M^p} = b^p \pi_0\right\} \le I\{b^p \pi_0 < 1\}\exp[-nK(b^p\pi_0, \pi_0)].$$

We complete the proof by inserting (A.10) into (A.9) and then inserting (A.9) into (A.8). □

PROOF OF PROPOSITION 6.3. The first inequality of (6.32) is that of (6.34). It follows from Lemma A.3 that $\sup_{\mu_p(G)\le C} R^*(G) - \mathcal{R}_n(\Theta_{p,n}(C))$ is bounded from above by a sum of three terms: two in Lemma A.3, and via the second inequality of (6.34), a third term bounded by

$$\text{(A.11)} \quad \sup_{\mu_p(G)\le C} R^*(G) - \sup_{\mu_p(G)\le C/b} R^*(G) \le (b^2-1)\sup_{\mu_p(G)\le C} R^*(G).$$

We choose $b$ and $\pi_0$ so that the three terms are of the same order.

Let $b^2 = 1 + \pi_0 \log_+(1/\pi_0)/\sup_{\mu_p(G)\le C} R^*(G)$. By Lemma A.3 and (A.11),

$$\text{(A.12)} \quad \sup_{\mu_p(G)\le C} R^*(G) - \mathcal{R}_n(\Theta_{p,n}(C))$$
$$\le (M'+1)\pi_0 \log_+(1/\pi_0) + \frac{I\{b^p\pi_0 < 1, b^2 < 2\}4C^2\pi_0^{-2/p}}{\exp[nK(b^p\pi_0, \pi_0)]}.$$

Since $K(p_1, p_2) \ge (p_1 - p_2)^2/(2p_1)$ for $p_2 < p_1 < 1$, for $1 < b^2 < 2$ and small $\pi_0 > 0$

$$K(b^p\pi_0, \pi_0) \ge \frac{(b^p\pi_0 - \pi_0)^2}{2(b^p\pi_0)} \ge \pi_0 b_0 p^2(b^2-1)^2 = \frac{b_0 p^2 \pi_0^3 \log_+^2(1/\pi_0)}{\{\sup_{\mu_p(G)\le C} R^*(G)\}^2},$$

where $b_0 \equiv \min[(b^p-1)^2/\{b^p p^2(b^2-1)^2\}: 1 \le b^2 \le 2, p > 0] > 0$. Thus, the second term in (A.12) is of the order $\pi_0$ for the choice of $\pi_0$ satisfying for certain $M'' < \infty$

$$b_0 n p^2 \pi_0^3 \log_+^2(1/\pi_0) \ge (M'')^2 \min\{1, C^{2p'}\{\log_+(1/C^{p'})\}^{2-p'}\}\log(C^2/\pi_0^{1+2/p}),$$

since $\sup_{\mu_p(G)\le C} R^*(G) \le M'' \min\{1, C^{p'}\{\log_+(1/C^{p'})\}^{1-p'/2}\}$ by (6.33). This holds with

$$\pi_0^3 \log_+^2(1/\pi_0) = 6(M'')^2 \min\{1, C^{2p'}\{\log_+(1/C^{p'})\}^{2-p'}\}\log_+(C^{p'})/(b_0 n p^2 p').$$

Hence, the conclusion holds, since $x^3 \log_+(1/x) = O(y)$ iff $x \log_+(1/x) = O(1)\Psi(y)$ for $x \wedge y > 0$. □

PROOF OF THEOREM 7.2. The proof of Theorem 3.2 provides an outline of the proof. The $o(\varepsilon^{-\eta})$ there is clearly bounded by a polynomial of $\log_+(1/\varepsilon)$. We omit the details, since the full proof of Theorem 7.2 can be found in Zhang (2000). □

PROOF OF THEOREM 7.3. The computation in the proof is similar to that in the proof of Theorem 7.2 and is provided in Zhang (2000). We again provide just an outline of the proof of (7.15) here. Let $L(\varepsilon)$ denote generic polynomials of $\log_+(1/\varepsilon)$. Let $B = B_{p,q}^{\alpha}(C)$ for fixed $(\alpha, p, q, C)$. By (7.12) and (6.31),

$$(\text{A.13}) \quad \sup_{\beta \in B} R^{(\varepsilon,*)}(\beta) - \mathcal{R}^{(\varepsilon)}(B) \leq L(\varepsilon)\varepsilon^2 \sum_{j=2}^{\infty} 2^j \tilde{r}_p^*(2^j, 2^{-j(\alpha+1/2)}/\varepsilon),$$

where $\tilde{r}_p^*(n, C) \equiv \min(1, C^{2p'/3})/n^{1/3}$. Splitting the sum in the right-hand side of (A.13) into two parts, for $2^{j(\alpha+1/2)}\varepsilon > 1$ and $\leq 1$, we find that the sum is of the order $\varepsilon^{-(2/3)/(\alpha+1/2)} = \varepsilon^{-1/(\alpha_3+1/2)}$. Thus, the left-hand side of (A.13) is bounded by $L(\varepsilon)\varepsilon^{2-1/(\alpha_3+1/2)} = L(\varepsilon)\varepsilon^{2\alpha_3/(\alpha_3+1/2)}$.

Now we prove (3.10) and (7.11). Let $j^* \geq 0$ satisfy $2^{j^*(\alpha+1/2)} \leq C/\varepsilon < 2^{(j^*+1)(\alpha+1/2)}$ and let $P$ be a probability measure under which $\beta_{j^*,k}$ are i.i.d. uniform variables in $[-\varepsilon, \varepsilon]$ and $\beta_{jk} = 0$ for $j \neq j^*$. By (3.9), $\|\beta\|_{p,q}^{\alpha} \leq 2^{j^*(\alpha+1/2-1/p)}\varepsilon 2^{j^*/p} \leq C$ almost surely under $P$, so that the minimax risk in $B_{p,q}^{\alpha}(C)$ is no smaller than the Bayes risk under $P$. With the scale change $\beta \to \beta/\varepsilon$, we find

$$\mathcal{R}^{(\varepsilon)}(B) \geq \inf_{\hat{\beta}} \int \left\{ \sum_{k=1}^{2^{j^*}} E_\beta^{(\varepsilon)}(\hat{\beta}_{j^*,k} - \beta_{j^*,k})^2 \right\} dP = 2^{j^*}\varepsilon^2 R^*(G_0)$$
$$\geq (C/\varepsilon)^{1/(\alpha+1/2)}\varepsilon^2 R^*(G_0)/2,$$

where $G_0$ is the uniform distribution in $[-1, 1]$ and $0 < R^*(G_0) < 1$ is the optimal Bayes risk in (6.1). This proves the lower bound in (3.10), and the lower bound, (7.10) and (7.15) imply (7.11). The upper bound in (3.10) follows from (7.10), (7.1), (7.4) and (6.33). □

PROOF OF THEOREM 8.1. Define $G_j(u) \equiv 2^{-j} \sum_k P\{\beta_{jk}/\varepsilon \leq u\}$, and define $G_j'$ and $G_j''$ in the same way for $\beta'$ and $\beta''$. Since $\overline{G}_j(u) \leq \overline{G}_j'(u/2) + \overline{G}_j''(u/2)$, by Proposition 6.1

$$r_0(2^j, G_j) \leq 2^p 3\, r_{p\wedge 2}(2^j, \mu_p(G_j')) + 4r_0(2^j, G_j'').$$

This splits the right-hand side of (8.1) into three sums. Since $\vec{\mu}_p(\beta') \in B_{p,q}^{\alpha}(C)$, (7.4) holds with $G_{[j]}^{(\varepsilon)} = G_j'$, so that $\varepsilon^2 \sum_j 2^j r_{p\wedge 2}(2^j, \mu_p(G_j')) = o(1)\varepsilon^{2\alpha/(\alpha+1/2)}$ as in the proof of Theorem 3.2. Moreover, since $\int_0^x \overline{G}_j''(\sqrt{u})\, du \leq x\kappa(\beta_{[j]}''/\varepsilon)$ for $x \geq 1$,

$$\sum_j 2^j r_0(2^j, G_j'') \leq \sum_j 2^j \log(2^j) \kappa(\beta_{[j]}''/\varepsilon)$$

$$\leq m^{(\varepsilon)} \sum_j \log(2^j)\left(1 \wedge \frac{M^{(\varepsilon)}}{\varepsilon^2 2^j}\right) \leq O(1)(\log \varepsilon)^2 m^{(\varepsilon)},$$

so that $\varepsilon^2 \sum_j 2^j r_0(2^j, G_j'')$ is also of order $o(1)\varepsilon^{2\alpha/(\alpha+1/2)}$. Thus, the right-hand side of (8.1) is uniformly $o(1)\varepsilon^{2\alpha/(\alpha+1/2)}$ over $\beta \in B^{(\varepsilon)}$. This proves (8.3).

It follows from (7.1) that $ER^{(\varepsilon,*)}(\beta) \leq \varepsilon^2 \sum_j 2^j R^*(EG_{[j]}^{(\varepsilon)}) = \varepsilon^2 \sum_j 2^j R^*(G_j)$. The total ideal risk for blocks with $2^j = o(1)\varepsilon^{-1/(\alpha+1/2)}$ is $o(1)\varepsilon^{2\alpha/(\alpha+1/2)}$. For blocks with $2^j \approx \varepsilon^{-1/(\alpha+1/2)}$, $R^*(G_j) = (1+o(1))R^*(G_j')$ by Proposition 6.2. For blocks with $\varepsilon^{-1/(\alpha+1/2)} = o(2^j)$, the total ideal risk is smaller than the optimal soft thresholding risk, which is $o(1)\varepsilon^{2\alpha/(\alpha+1/2)}$ as in the proof of (8.3). Thus, (8.4) holds. We omit certain details. □

PROOF OF THEOREMS 4.1 AND 9.1. We first prove Theorem 9.1(i). It follows from the proof of Lemma 7 in Brown, Cai, Low and Zhang (2002) that

$$E_f(\tilde{\beta}_{j,k} - \beta_{j,k})^2 \leq 4\frac{2^{j\vee 0}}{N}\int(f - \bar{f}_{j\vee 0})^2 \mathbb{1}_{j,k} - 3\frac{2^{j\vee 0}}{N}\int(f - \bar{f}_{j+1})^2 \mathbb{1}_{j,k}$$

(A.14)

$$= 4\frac{2^j}{N}\beta_{j,k}^2 I\{j \geq 0\} + \frac{2^{j\vee 0}}{N}\sum_{\ell=j+1}^{\infty}\sum_{m=1}^{2^\ell}\beta_{\ell,m}^2 \mathbb{1}_{j,k}(m/2^\ell),$$

since $(\bar{f}_{\ell+1} - \bar{f}_\ell)\mathbb{1}_{\ell,m} = \beta_{\ell,m}\phi_{\ell,m}$ and $\int |\bar{f}_{\ell+1} - \bar{f}_\ell|^2 \mathbb{1}_{\ell,m} = \beta_{\ell,m}^2$. Thus, for $p' \equiv p \wedge 2$

(A.15)
$$\sum_{k=1}^{2^{j\vee 1}} E_f|\tilde{\beta}_{j,k} - \beta_{j,k}|^{p'}$$

$$\leq \frac{2^{(j\vee 0)p'/2}}{N^{p'/2}}\left\{2^{p'}\sum_{k=1}^{2^{j\vee 1}}|\beta_{j,k}|^{p'}I\{j \geq 0\} + \sum_{\ell=j+1}^{\infty}\sum_{m=1}^{2^\ell}|\beta_{\ell,m}|^{p'}\right\}.$$

Since $\sum_{m=1}^{2^{j\vee 1}} |\beta_{j,m}|^{p'} \leq 2^{-(j\vee 0)p'(\alpha+1/2-1/p')}(\|\beta\|_{p',q}^\alpha)^{p'}$ and $\|\beta\|_{p'q}^\alpha \leq \|\beta\|_{p,q}^\alpha$, by (A.15)

$$2^{(j\vee 0)(\alpha+1/2-1/p')}\left(E_f\sum_{k=1}^{2^{j\vee 1}}|\tilde{\beta}_{j,k} - \beta_{j,k}|^{p'}\right)^{1/p'} \leq M_{\alpha,p'}\|\beta\|_{p,q}^\alpha\sqrt{2^j/N},$$

so that (9.1) holds. Furthermore, (A.15) with $p = 2$ implies that $E_f\|\Pi_J\tilde{\beta} - \beta\|_2^2$ is bounded by

$$\left\{E_f|\tilde{\beta}_{-1,1} - \beta_{-1,1}|^2 + E_f\sum_{j=0}^{J}\sum_{k=1}^{2^j}|\tilde{\beta}_{j,k} - \beta_{j,k}|^2 + \sum_{j=J+1}^{\infty}\sum_{k=1}^{2^j}\beta_{j,k}^2\right\}$$

$$\leq \sum_{\ell=0}^{\infty} \Big(\frac{2^{\ell \wedge (J+1)}}{N} + \frac{2^{2+\ell}}{N} I\{\ell \leq J\} + I\{\ell > J\}\Big) \Big(\sum_{m=1}^{2^\ell} |\beta_{\ell,m}|^{p'}\Big)^{2/p'}$$

$$\leq (\|\beta\|_{p,q}^{\alpha})^2 \bigg[\frac{5}{N} \sum_{\ell=0}^{J} 2^{\ell\{1-2(\alpha+1/2-1/p')\}}$$

$$+ \Big(\frac{2^{J+1}}{N} + 1\Big) \sum_{\ell=J+1}^{\infty} 2^{-2\ell(\alpha+1/2-1/p')}\bigg]$$

$$\leq M_{\alpha,p'}(\|\beta\|_{p,q}^{\alpha})^2 \Big\{\frac{J}{N} I\{\alpha p' = 1\} + \frac{1}{N} + \Big(\frac{2^J}{N}+1\Big) 2^{-2J(\alpha+1/2-1/p')}\Big\}.$$

This implies (9.2) and completes the proof of Theorem 9.1(i).

Now we prove that for $\varepsilon = \sigma/\sqrt{N}$ and the $\zeta_N$ in Theorem 4.1

(A.16) $$\sup_{f \in \mathcal{F}} E_f \|\hat{f}_N - f\|^2 \leq (1+\zeta_N) \mathcal{R}^{(\varepsilon)}(B_{p,q}^{\alpha}(C)).$$

Define $\widetilde{G}'_j(u) = n_j^{-1} \sum_k \delta_{j,k} I\{\tilde{\beta}_{j,k} \leq u\}$ and $\widetilde{G}_j(u) = 2^{-j} \sum_k P\{\tilde{\beta}_{j,k} \leq u\}$ with the $\tilde{\beta}_{j,k}$ in (4.7). Let $\tilde{y}_{j,k} \equiv \delta_{j,k} y_{j,k} + (1-\delta_{j,k}) N(0, \varepsilon^2)$. By Theorem 3.1 and (3.3)

$$E_f \sum_{j=-1}^{J} \sum_k (\hat{\beta}_{j,k} - \tilde{\beta}_{j,k})^2$$

$$\leq E_f \sum_{j=-1}^{J} \inf_{\{t_j\}} \sum_k \delta_{j,k} (t_j(y_{j,k}) - \tilde{\beta}_{j,k})^2 + \varepsilon^2 \sum_{j=-1}^{J} E_f n_j r_{p'}(n_j, \mu_{p'}(\widetilde{G}'_j))$$

$$\leq E_f \sum_{j=-1}^{J} \inf_{\{t_j\}} \sum_k (t_j(\tilde{y}_{jk}) - \tilde{\beta}_{j,k})^2 + \varepsilon^2 \sum_{j=-1}^{J} 2^j r_{p'}(2^j, \mu_{p'}(\widetilde{G}_j)).$$

Since $\tilde{y}_{j,k} \sim N(\tilde{\beta}_{j,k}, \varepsilon^2)$ given $t_i$, a slight modification of the proof of Theorem 8.1 implies that the right-hand side above is bounded by $(1+\zeta_N)\mathcal{R}^{(\varepsilon)}(B_{p,q}^{\alpha}(C))$, in view of Theorem 9.1(i). This and (9.2) imply (A.16) with the $\zeta_N$ in Theorem 4.1 for the choices of $J = J_N$ in Theorem 4.1.

It remains to show that for $\alpha + 1/2 - 1/p > \alpha/(\alpha+1/2)$ and $\tilde{f}$ based on (4.1)

(A.17) $$\inf_{\tilde{f}} \sup_{f \in \mathcal{F}} E_f \|\tilde{f} - f\|^2 \geq (1+o(1))\mathcal{R}^{(\varepsilon)}(B_{p,q}^{\alpha}(C)).$$

Let $T_N \equiv T_{N,\{t_i\}}$ be the randomized mappings $\{Y_i, t_i, i \leq N\} \to \{\tilde{y}_{j,k}, \delta_{j,k}\}$. Brown, Cai, Low and Zhang (2002) proved that due to the orthonormality of the mappings $T_N$ given $\{t_i\}$, the inverse mappings of $T_N$ provide $\{\tilde{y}_{j,k}, \delta_{j,k}\} \to \{Y'_i, t_i, i \leq N\}$ satisfying (4.1) with regression functions $f'(t)$ such that (A.14) holds with $\tilde{\beta}_{j,k} = \int f' \phi_{jk}$. This yields (A.17) by repeating the proof of (A.16). $\square$

DEPARTMENT OF STATISTICS
RUTGERS UNIVERSITY
HILL CENTER, BUSCH CAMPUS
PISCATAWAY, NEW JERSEY 08854
USA
E-MAIL: czhang@stat.rutgers.edu